\documentclass[12pt]{elsarticle}

\usepackage[latin1]{inputenc}
\usepackage{mathstyle00}
\usepackage{bm}
\usepackage[margin=0.98in]{geometry}
\begin{document}
\begin{frontmatter}

\title{A class of robust numerical methods for solving dynamical systems with multiple time scales}
    
  \author[cms]{Thomas Y. Hou}
  \ead{hou@cms.caltech.edu}			 
  \author[hku]{Zhongjian Wang}
  \ead{ariswang@connect.hku.hk}
  \author[hku]{Zhiwen Zhang\corref{cor1}}
  \ead{zhangzw@hku.hk}
            
            \address[cms]{Computing and Mathematical Sciences, California Institute of Technology, Pasadena, CA 91125, USA. }
            \address[hku]{Department of Mathematics, The University of Hong Kong, Pokfulam Road, Hong Kong SAR, China.}
            \cortext[cor1]{Corresponding author}
            
    \begin{abstract}
    \noindent In this paper, we develop a class of robust numerical methods for 
    solving dynamical systems with multiple time scales. We first represent the solution of a multiscale dynamical system as a transformation of a slowly varying solution. Then, under the scale separation assumption, we provide a systematic way to construct the transformation map and derive the dynamic equation for the slowly varying solution. We also provide the convergence analysis of the proposed method.
    Finally, we present several numerical examples, including ODE system with three and four separated time scales to demonstrate the accuracy and efficiency of the proposed method. Numerical results verify that our method is robust in solving ODE systems with multiple time scale, where the time step does not depend on the multiscale parameters.

\end{abstract}
\end{frontmatter}

\section{Introduction} \label{sec:introduction}
\noindent
Dynamical systems with sub-processes evolving on many different time scales are ubiquitous in  applications: chemical reactions, electro-optical and neuro-biological systems, to name just  a few \cite{tuckerman1992reversible,jones2012multiple}. The multiple time scales in the dynamical systems pose a major problem in numerical simulations because one needs to choose small time steps for stable integration of the fast motions in the systems, which leads to large numbers of time steps required for the observation of the slow degrees of freedom and thus requires tremendous computational resources. 
Interested readers are referred to \cite{hairer2006geometric,kuehn2015multiple} and the references therein for a detailed review.  

The objective of this paper is to develop a new method to solve dynamical systems with multiple time scales; see Section \ref{sec:NumSchemeMultiscaleDS} for the precise definition of the problems. The main idea of our method is to formally represent the solution of the multiscale dynamical system as a transformation of a slowly varying solution; see Eq.\eqref{xmapy}. By dealing with the multiscale information in a dimension-by-dimension fashion, we propose a systematic way to construct a set of cumulative composition maps that capture the complicated dynamics of the problem. 
Based on the scale separation assumption, we successfully derive the dynamic equation for the slowly 
varying solution (i.e., Eq.\eqref{eqn:non-stiff-nscale-1D}) and prove that the dynamic equation for the slowly varying solution is non-stiff; see Theorem \ref{thm:nonstiff}. Thus, we can use conventional numerical methods to compute it, where the time step is independent of the multiscale parameters in the dynamical system. In addition, we analyze the error between the numerical solution obtained from our method and the exact solution in Theorem \ref{thm:mainerror}. Finally, we carry out several numerical experiments to demonstrate the accuracy and efficiency of the proposed method. 

As we will demonstrate in Section \ref{sec:numericalresults}, the proposed method can offer accurate numerical solutions to multiscale ODE systems with considerable computational savings over traditional methods, especially when the multiscale parameters are small. Numerical results (see Fig.\ref{fig:x1x5_avg}) show that the dynamic equation for the slowly varying solution based on the cumulative composition maps indeed capture the averaged behaviors of the solution well. While a simple averaging treatment of the original multiscale ODE systems leads to wrong results. As an analogy to this interesting finding, in the homogenization for elliptic PDEs with multiscale coefficients,  a simple average of the coefficient gives a wrong result, where one needs to solve a cell problem to obtain the correct homogenization coefficient \cite{EfendievHou:09}. 

Our method is inspired by the recent development in designing uniformly accurate numerical schemes for highly oscillatory evolution equations \cite{chartier2015uniformly,chartier2017new}, where two-scale 
problems were solved. In \cite{chartier2015uniformly,chartier2017new}, the authors separate the two time scales into two independent variables and embed the solution of the two-scale problem into a two-variable function. Then, they derive formulations of the evolution equations for the two-variable function and prove that under certain conditions the evolution equations are solvable and non-stiff. 

The novelty of our paper is that we provide a systematical way to construct a set of cumulative composition maps that allow us to correctly upscale the complicated dynamics of the problem. Notice from Eq.\eqref{defi:phik1D} that each map $\Phi^k$ is a perturbation of the identity operator. However, 
a cumulative composition of those simple maps \eqref{cumcompmaps} can provide an accurate approximation of the complicated dynamics of the problem.  In addition, we provide a rigorous convergence analysis for the proposed method and verify the statement through numerical experiments.  
     
Before we end this section, we give a short review of several existing methods for solving two-scale 
problems. When slow variables can be identified, effective equations can be obtained by averaging the instantaneous drift driving those slow variables. Two classes of numerical methods have been developed based on this observation: the equation-free method \cite{kevrekidis2003equation} and heterogeneous multiscale method (HMM) \cite{abdulle2012heterogeneous}. Later on, a new class of integrators for stiff ODEs as well as SDEs 
were developed \cite{tao2010nonintrusive}, which are based on the averaging of the instantaneous flow of 
the hidden slow and fast variables simultaneously. Therefore, the hidden slow variables do not need to be explicitly identified.  In this paper, however, we will consider problems parameterized by multiple time scales.  In addition, we aim to design numerical schemes that solve the multiscale dynamical problems for a wide range of multiscale parameter values with uniform accuracy.

The rest of the paper is organized as follows. In Section 2, we will derive our numerical method  
for solving multiscale dynamical systems and discuss its detailed implementation. In Section 3, we provide the convergence analysis for the proposed method.  In Section 4, we present numerical results to demonstrate the accuracy and efficiency of our method.  Concluding remarks are made in Section 5.
 
\section{Numerical methods for solving multiscale dynamical systems}\label{sec:NumSchemeMultiscaleDS}
\noindent	
In this section, we will develop numerical methods to solve dynamical systems with multiple time scales. Specifically, we consider the following first-order ordinary differential equation (ODE) system to illustrate the main idea,
\begin{equation}
\dot{x}=f^{\bm{\epsilon}}(t,x),\quad x(0)=x_0, \quad t\in[0,T],
\label{ODEsystem-multiscale}
\end{equation} 
where $x(t)\in \mathbb{R}^{d}$ is the solution vector, $x_0$ is the initial value, and $f^{\bm{\epsilon}}(t,x)$ is a function vector field. Here $\bm{\epsilon}=(\epsilon_1,...,\epsilon_n)$ is a set of parameters, which are used to characterize different time scales in the ODE system \eqref{ODEsystem-multiscale}. { When the parameters satisfy
\begin{equation}
0<\epsilon_n\ll\epsilon_{n-1}\ll\cdot\cdot\cdots\ll\epsilon_1\ll1, 
\label{WellSeparated-multiscale}
\end{equation}
we say that the multiscale time scales are well-separated.}  Given the multiscale parameters, we denote 
\begin{equation}
f^{\bm{\epsilon}}(t,x)\equiv f(\frac{t}{\epsilon_1},\frac{t}{\epsilon_2},\cdots,\frac{t}{\epsilon_n},x), 
\label{f-multiscales}
\end{equation}
where $f:\mathbb{R}^{n+d}\to \mathbb{R}^{d}$ is a function vector field.  We assume that the first-order derivatives of $f$ are bounded, which is sufficient to guarantee the existence and uniqueness of solutions of the ODE system \eqref{ODEsystem-multiscale} \cite{iserles2009first,perko2013differential}. Let $t_i=\frac{t}{\epsilon_i}$, $i=1,...,n$. We denote $f^{\bm{\epsilon}}(t,x)=f(t_1,t_2,\cdots,t_n,x)$. Moreover, we assume that $f$ is periodic with respect to its first $n$ coordinates, i.e., $t_i$, $i=1,...,n$. Without lost of generalities, all the periods are assumed to be $1$.  

\subsection{Decomposition of the multiscale function $f$}\label{sec:CompositionMaps}
\noindent
We iteratively define the averaged functions to resolve finer scale fluctuations on coarser scales in a dimension-by-dimension fashion. We first start from the coordinate $t_n$ corresponding to the smallest-scale and define the mean function
\begin{equation}
\bar{f}^n(t_1,t_2,\cdots,t_{n-1},x)=\int_{0}^{1}f(t_1,t_2,\cdots,t_{n-1},s,x)ds,
\label{mean-fn}
\end{equation} 
and the fluctuation function 
\begin{equation}
f^{n}(t_1,t_2,\cdots,t_{n},x)=f(t_1,t_2,\cdots,t_{n},x)-\bar{f}^n(t_1,t_2,\cdots,t_{n-1},x),
\label{fluc-fn}
\end{equation}
where the integration and subtraction are done in a component-wise fashion. Thus, $\bar{f}^n$ 
and $f^{n}$ are $d$-dimensional vector functions. 

Then, based on the function $\bar{f}^n(t_1,t_2,\cdots,t_{n-1},x)$, we define the mean function and  fluctuation function corresponding to the second smallest-scale as follows 
\begin{align}
\bar{f}^{n-1}(t_1,t_2,\cdots,t_{n-2},x)&=\int_{0}^{1}\bar{f}^n(t_1,t_2,\cdots,t_{n-2},s,x)ds,\label{mean-fn1}\\
f^{n-1}(t_1,t_2,\cdots,t_{n-1},x)&=\bar{f}^n(t_1,t_2,\cdots,t_{n-1},x)-\bar{f}^{n-1}(t_1,t_2,\cdots,t_{n-2},x).
\label{fluc-fn1}
\end{align}
We continue this strategy and define mean functions and fluctuation functions corresponding to different time scales recursively. For instance, given the mean function $\bar{f}^{n-k+1}$, we define the mean function and fluctuation function corresponding to a coarser-scale as follows
\begin{align}
\bar{f}^{n-k}(t_1,t_2,\cdots,t_{n-k-1},x)&=\int_{0}^{1}\bar{f}^{n-k+1}(t_1,t_2,\cdots,t_{n-k-1},s,x)ds,
\label{mean-fnk}\\
f^{n-k}(t_1,t_2,\cdots,t_{n-k},x)&=\bar{f}^{n-k+1}(t_1,t_2,\cdots,t_{n-k},x)
-\bar{f}^{n-k}(t_1,t_2,\cdots,t_{n-k-1},x).
\label{fluc-fnk}
\end{align}
Finally, we define the mean function and fluctuation function corresponding to the largest-scale as follows
\begin{align}
\bar{f}^{1}(x)&=\int_{0}^{1}\bar{f}^{2}(s,x)ds, \label{mean-f1}\\
f^{1}(t_1,x)&=\bar{f}^{2}(t_1,x)-\bar{f}^{1}(x).\label{fluc-f1}
\end{align}

	The above recursive formulations \eqref{mean-fn}-\eqref{fluc-f1} naturally lead to a decomposition of the multiscale function $f(t_1,t_2,\cdots,t_n,x)$ into 
	\begin{equation}
	f=\sum_{k=1}^{n}f^k+\bar{f}^1.
	\label{nest-decomp-f}
	\end{equation}
	According to the definition, for each $k$ we have that
	\begin{equation}\label{fkmeanzero}
	\int_{0}^1f^{k}(t_1,\cdots,t_{k-1},s,x)ds=0, \quad \forall~t_i\in[0,1],~i=1,...,k-1,~\text{and}~ x.
	\end{equation}
	\subsection{Derivation of the dynamic equation for the slowly varying solution}\label{sec:DynamicSlowEq}
	\noindent
	We will construct a family of maps $\Phi^k=\Phi^k_{t_1,t_2,\cdots,t_k}:\mathbb{R}^d\to\mathbb{R}^d$, $k=1,...,n$, that allow us to represent the solution $x(t)$ of the ODE system \eqref{ODEsystem-multiscale} as a transformation of a slowly varying solution $y(t)$, i.e.,   
	\begin{equation}
	x(t)=\Phi^n\circ\Phi^{n-1}\circ\cdots\circ\Phi^1(y(t)).
	\label{xmapy}
	\end{equation} 
	We assume that each map $\Phi^k$ is periodic with respect to $t_k$ and becomes the identical map $I_d$ when $t_k=0$. To simplify the notation,  we define a family of cumulative composition maps as follows 
	\begin{equation}
	\bar{\Phi}^k=\Phi^k\circ\Phi^{k-1}\circ\cdots\circ\Phi^1, \quad k=1,...,n.
	\label{cumcompmaps}
	\end{equation} 
	
	The rationale behind the representation \eqref{xmapy} is that the complicated dynamics of the ODE system \eqref{ODEsystem-multiscale} (e.g., highly oscillatory solutions) is captured by the map $\bar{\Phi}^n$, and thus the evolution of the solution $y(t)$ is smooth. Therefore, we can compute the solution $x(t)$ (through solving $y(t)$) by using numerical methods with relatively large time steps (independent of the multiscale parameters). To achieve this goal, we need to find a constructive way to obtain the map $\bar{\Phi}^n$ and to identify the dynamical equation for the solution $y(t)$. 
	
	Substituting the representation \eqref{xmapy} into the original problem \eqref{ODEsystem-multiscale}, we know that the solution $y(t)$ formally satisfies the following equation 
	\begin{equation}
	\partial_t\bar{\Phi}^n(y)+\partial_{x}\bar{\Phi}^n(y)\dot{y}=f^\epsilon, 
	\label{eqn:governing-nph-1D}
	\end{equation}
	where $\partial_{x}\bar{\Phi}^n(y)$ is the Jacobian matrix. 
	
	It will be complicated if we directly compute $\partial_{x}\bar{\Phi}^n(y)$ using the chain rule. We shall adopt some approximate method to address this difficulty.  Notice that the scale separation assumption on $f^\epsilon$ (see Eq.\eqref{f-multiscales}) indicates that the local fluctuation of $f^\epsilon$ at $\mathcal{O}(\epsilon_n)$ scale can be resolved by the same scale part of $\partial_t\bar{\Phi}^n$. More precisely, by the chain rule, we get 
	\begin{align}\label{eqn:Order-n-expansion}
	\partial_t\bar{\Phi}^n=\frac{1}{\epsilon_n}\partial_{t_n}\Phi^n(\bar{\Phi}^{n-1})+\sum_{i=1}^{n-1}\frac{1}{\epsilon_{i}}\partial_{t_i}\Phi^n(\bar{\Phi}^{n-1})+\partial_x\Phi^n\partial_{t}\bar{\Phi}^{n-1},
	\end{align}
	where the last two terms are independent of $t_n$. Compared with the decomposition of the multiscale function $f$ in \eqref{nest-decomp-f}, we can set,
	\begin{equation}
	\frac{1}{\epsilon_n}\partial_{t_n}\Phi^n(\bar{\Phi}^{n-1})=f^n.
	\end{equation}
	Since $\Phi^n$ is an identical map when $t_n=0$, we can explicitly get,
	\begin{equation}
	\Phi^n_{t_1,t_2,\cdots,t_n}(x)=x+\epsilon_ng^n(t_1,t_2,\cdots,t_n,x),\quad \forall x,
	\label{Phin}
	\end{equation}
	where \begin{equation}
	g^{n}(t_1,t_2,\cdots,t_n,x)=\int_{0}^{t_n}f^{n}(t_1,t_2,\cdots,t_{n-1},s,x)ds.
	\label{gintf}
	\end{equation} 
	From \eqref{Phin}, we can see that $\Phi^n=I_d+\mathcal{O}(\epsilon_{n})$, which is an $\mathcal{O}(\epsilon_{n})$ order perturbation of the identity operator. When the time scales are well-separated, we can dismiss the second term in Eq.\eqref{eqn:Order-n-expansion} since 
	\begin{equation}
	\sum_{i=1}^{n-1}\frac{1}{\epsilon_{i}}\partial_{t_i}\Phi^n(\bar{\Phi}^{n-1})=\sum_{i=1}^{n-1}\frac{\epsilon_{n}}{\epsilon_{i}}\partial_{t_i}g^n(\bar{\Phi}^{n-1})=\mathcal{O}(\frac{\epsilon_{n}}{\epsilon_{n-1}}).
	\end{equation}  Let us continue our derivation inductively with $k=n-1,\cdots,1$. We consider the fluctuation within period $\mathcal{O}(\epsilon_k)$, and have the following observation, 
	\begin{equation}
	\partial_x(\Phi^n\circ\Phi^{n-1}\circ\cdots\circ\Phi^{k+1})(\frac{1}{\epsilon_{k}}\partial_{t_k}\Phi^k)=f^k.
	\end{equation}
	Due to the scale separation structure of $\Phi^n\circ\Phi^{n-1}\circ\cdots\circ\Phi^{k+1}$, we obtain
	\begin{equation}
	\partial_x(\Phi^n\circ\Phi^{n-1}\circ\cdots\circ\Phi^{k+1})=I_d+\mathcal{O}(\epsilon_{k+1}), 
	\end{equation}
	which is an $\mathcal{O}(\epsilon_{k+1})$ order perturbation of the identity operator. Now we arrive at, 
	\begin{equation}
	\frac{1}{\epsilon_k}\partial_{t_k}\Phi^k(\bar{\Phi}^{k-1})=f^k.
	\end{equation}
	Again, using the condition that $\Phi^k$ is an identical map when $t_k=0$, we get  
	\begin{equation}\label{defi:phik1D}
	\Phi^k_{t_1,t_2,\cdots,t_k}(x)=x+\epsilon_kg^k(t_1,\cdots,t_{k-1},t_k,x),
	\end{equation}
	where \begin{equation}
	g^{k}(t_1,\cdots,t_{k-1},t_k,x)=\int_{0}^{t_k}f^{k}(t_1,\cdots,t_{k-1},s,x)ds, \quad k=1,...,n.
	\label{gkintfk}
	\end{equation}
	In the above derivation, we have used the condition that  
	$\sum_{i=1}^{k-1}\frac{1}{\epsilon_{i}}\partial_{t_i}\Phi^k(\bar{\Phi}^{k-1})=\mathcal{O}(\frac{\epsilon_{k}}{\epsilon_{k-1}})$.
	
	After we obtain the explicit formulations for $\Phi^k$, $k=1,...,n$ and their derivatives, we are in the position to derive the dynamic equation for the solution $y(t)$.  
	According to \eqref{eqn:governing-nph-1D}, we obtain a nested equation 
	\begin{align}
	\partial_{t}\Phi^n+\partial_x\Phi^n\bigg(\partial_{t}\Phi^{n-1}+\partial_x\Phi^{n-1}\Big(\partial_{t}\Phi^{n-2}+\cdots(\partial_{t}\Phi^{1}+\partial_x\Phi^{1}\dot{y}(t))\cdots\Big)\bigg)=f^\epsilon.
	\label{master-eq1}
	\end{align}
	From the definition of $\Phi^k$ (see Eq.\eqref{defi:phik1D}), we compute the derivative of $\Phi^k$ with respect to time $t$ and get,
	\begin{align}
	\partial_{t}\Phi^k=\epsilon_{k}\sum_{r=1}^{k}\frac{1}{\epsilon_r}\partial_{t_r}g^k, \quad k=1,...,n.
	\end{align}
	The scale separation assumption on the multiscale parameters implies that   $\partial_{t}\Phi^k\approxeq\partial_{t_k}g^k=\frac{1}{\epsilon_k}\partial_{t_k}\Phi^k$, which 
	allows us to simplify Eq.\eqref{master-eq1} into the following form
	\begin{align}
	\frac{1}{\epsilon_n}\partial_{t_n}\Phi^n+&\partial_x\Phi^n\bigg(\frac{1}{\epsilon_{n-1}}\partial_{t_{n-1}}\Phi^{n-1}+\partial_x\Phi^{n-1}\cdot\nonumber\\&\Big(\frac{1}{\epsilon_{n-2}}\partial_{t_{n-2}}\Phi^{n-2}+\cdots(\frac{1}{\epsilon_{1}}\partial_{t_{1}}\Phi^{1}+\partial_x\Phi^{1}\dot{\tilde{y}}(t))\cdots\Big)\bigg)=f^\epsilon.
	\label{master-eq2}
	\end{align}
	Here $\tilde{y}(t)$ is an approximation of $y(t)$ since Eq.\eqref{master-eq2} is an approximation 
	of the original Eq.\eqref{master-eq1} based on the scale separation assumption on the multiscale parameters. Finally, from Eq.\eqref{master-eq2} we can get the dynamic equation for $\tilde{y}(t)$, i.e.,  
	\begin{align}\label{eqn:non-stiff-nscale-1D}
	\dot{\tilde{y}}=&(\partial_x\Phi^{1})^{-1}\bigg(\cdots(\partial_x\Phi^{n-1})^{-1}\nonumber\\&\Big((\partial_x\Phi^{n})^{-1}(f^\epsilon-\frac{1}{\epsilon_{n}}\partial_{t_{n}}\Phi^{n})-\frac{1}{\epsilon_{n-1}}\partial_{t_{n-1}}\Phi^{n-1}\Big)-\cdots-\frac{1}{\epsilon_{1}}\partial_{t_{1}}\Phi^{1}\bigg)=:F(t,x).
	\end{align}	
{ 
From the above derivation, one can see that the existence of the matrices $(\partial_x\Phi^{k})^{-1}$, $k=1,...,n$ in \eqref{eqn:non-stiff-nscale-1D} is essential in establishing the consistency of our method. In Section \ref{sec:convergence-analysis}, we will prove that the invertibility is guaranteed in the case when $\epsilon_{k}$ are sufficient small; see \eqref{Phikinv}}.	In addition, we will prove that the ODE system \eqref{eqn:non-stiff-nscale-1D} is non-stiff, which will be useful for the design of uniformly accurate numerical schemes, i.e., the time step in the numerical schemes does not depend on the multiscale parameters. When we obtain the solution  $\tilde{y}(t)$ of the ODE system \eqref{eqn:non-stiff-nscale-1D}, we can recover the solution of the ODE system \eqref{ODEsystem-multiscale} through the transform $\bar{\Phi}^n$ defined in \eqref{xmapy}\eqref{cumcompmaps}. The error estimate of our method will be presented later.  
\begin{remark}
From the explicit formulations for $\Phi^{k}$, $k=1,...,n$, we know the solution of the original 
ODE system can be rewritten as the following form
\begin{equation}
x(t)=\big(I_d+\mathcal{O}(\epsilon_{n})\big)\circ\big(I_d+\mathcal{O}(\epsilon_{n-1})\big) \circ\cdots\circ\big(I_d+\mathcal{O}(\epsilon_{1})\big)(y(t)).
\label{cumcompmapexplain}
\end{equation} 
Eq.\eqref{cumcompmapexplain} clearly reveals the structure of the transformation map in our method. 
One can see that the transformation map is a composition of simple maps, where each of them is a perturbation of identity. Interestingly, similar ideas appear in deep neural network research; see e.g. \cite{bartlett2018representing,shen2019nonlinear}, where approximations of functions via compositions of near-identity functions { have been used intensively and are the key to the amazing expressibility power of a deep neuron network.}  
\end{remark}


\begin{remark}\label{quasi-periodic}
Our method can be extended to solve an ODE system \eqref{ODEsystem-multiscale}, where $f^{\bm{\epsilon}}$ is a quasi-periodic function.
Assume that $f^{\bm{\epsilon}}(t,x)$ in \eqref{f-multiscales} has the form
\begin{align}
f^{\bm{\epsilon}} (t,x) = f(\frac{a_1(t)}{\epsilon_1}, \frac{a_2(t)}{\epsilon_2}, ..., \frac{a_n(t)}{\epsilon_n}, x), \label{quasi-periodic}
\end{align} 
where $a_k(t)$, $k=1,...,n$ are some invertible functions in $C^2$ such that 
$0< c_0 \le \| \frac{d}{dt} a_k(t)\|_{\infty} \le c_1 < \infty$. 
Then the main results stated in this section 
for periodic functions still hold by using the same definition of the mean function defined in \eqref{mean-fn} without any prior knowledge of $a_k(t)$, $k=1,...,n$. The reason is that for any smooth function $h(x,y)$ that is periodic in $y$ with period $1$, one can easily show that (see \cite{engquist1989particle} for an elementary proof) 
\begin{align}
\Big| \int_a^b h \big( \frac{a_{k-1}(t)}{\epsilon_{k-1}}, \frac{a_{k}(t)}{\epsilon_{k}}\big) dt - 
\int_a^b\big(\int_0^1 h ( \frac{a_{k-1}(t)}{\epsilon_{k-1}}, y) dy\big)dt \Big| \le C \frac{\epsilon_k}{\epsilon_{k-1}},
\end{align}
by using a change of variable from $t$ to $s=a_k(t)$ and the fact that the Jacobian 
$J(s) = (\frac{d}{dt} a_k)^{-1}$ is a smooth function of $s$.

\end{remark}

\begin{remark}\label{dealwithnonscale} 
	For a general ODE system $\dot{x}=f(t,x)$, where $f(t,x)$ does not have an explicit form of multiscale separation parametrization, we may reparameterize and approximate $f(t,x)$ by a formal multiscale velocity field $f^{\bm{\epsilon}}(t,x)$. For instance, we may reparameterize $f(t,x)$ into a formal two-scale structure through Fourier transform; see \cite{hou2008multiscale}. The limitation is that we have to compute Fourier transform of $f(t,x)$ with respect to $t$, for each fixed $x$, which involves a certain amount of computation. To develop a fast solver to address this issue will be our future work.     
\end{remark}

\subsection{Construction of the numerical schemes}\label{sec:Implentmentation}
\noindent
In this section, we construct efficient numerical schemes to solve Eq.\eqref{eqn:non-stiff-nscale-1D}
that are uniformly accurate with respect to $\bm{\epsilon}=(\epsilon_1,...,\epsilon_n)$. We first discuss how to accurately and efficiently compute the maps $\Phi^k$, $k=1,...,n$ defined in Eqns.\eqref{Phin} and \eqref{defi:phik1D}. We observe that the maps $\Phi^k$, $k=1,...,n$ are explicitly defined. Therefore 
we can use an explicit numerical scheme to approximate them. However, such an explicit implementation 
is not desirable because it may destroy the structures (e.g., Hamiltonian structure) of the original problem \eqref{ODEsystem-multiscale}.

Alternatively, we adopt an implicit midpoint scheme to approximate $\Phi^k$, i.e.,
\begin{equation}\label{midpoint-Phik}
\Phi^k_{t_1,t_2,\cdots,t_k}(x)=x+\epsilon_kg^k(t_1,t_2,\cdots,t_k,\frac{x+\Phi^k_{t_1,t_2,\cdots,t_k}(x)}{2}),\quad k=1,...,n. 
\end{equation}
The scheme \eqref{midpoint-Phik} still provides an $\mathcal{O}(\epsilon_k)$ approximation of Eq.\eqref{defi:phik1D}. 
In practice,  $\Phi^k$ in \eqref{midpoint-Phik} can be computed by the fixed point iteration. 
In addition, the derivatives of $\Phi^k$ with respect to $t_k$ or $x$ involved in Eq.\eqref{eqn:non-stiff-nscale-1D} can be computed by the fixed point iteration based on the following identities, 
\begin{align}\label{eqn:iterative-nscale}
\partial_{t_k}\Phi^k(x)&=\epsilon_{k}f^k\big(\frac{x+\Phi^k(x)}{2}\big)+\frac{\epsilon_{k}}{2}\partial_{x}g^k\big(\frac{x+\Phi^k(x)}{2}\big)\partial_{t_k}\Phi^k(x),\\
\big(\partial_x\Phi^k(x)\big)^{-1}K&=K-\frac{\epsilon_{k}}{2}\partial_{x}g^k\big(\frac{x+\Phi^k(x)}{2}\big)\big(K+(\partial_x\Phi^k(x))^{-1}K\big)\label{eqn:iterative-nscale2},
\end{align}
where $K$ is a {  $d$-dimensional} column vector. The formulae in \eqref{midpoint-Phik}-\eqref{eqn:iterative-nscale2} suggests an iterative scheme to
calculate all the quantities that are needed to compute \eqref{eqn:non-stiff-nscale-1D} and $\bar{\Phi}^n$. Finally, we obtain an efficient numerical scheme to solve Eq.\eqref{eqn:non-stiff-nscale-1D} at any time $t$ and value $x$. 

The detailed implementation of the proposed numerical scheme is listed in Algorithm \ref{alg:nscale}, in which we introduce several variables to simply the notations. Specifically, we have $P_k= \Phi^k\circ\Phi^{k-1}\circ\cdots\circ\Phi^1(y)$, $T_k=\frac{1}{\epsilon_k}\partial_{t_k}\Phi^k\circ\Phi^{k-1}\circ\cdots\circ\Phi^1(y)$, $k=1,...,n$, and $D_1= \dot{y}$.

\begin{algorithm}[ht]
	\caption{A fixed point iteration method to compute the ODE with $n$ time-scales.}
	\label{alg:nscale}
	\begin{algorithmic}[1]
		\STATE {Set $i=0$, $P_1^{[0]}=P_2^{[0]}=\cdots=P_n^{[0]}=y$} 
		\REPEAT
		\STATE $P_0^{[i]}=y$
		\FOR{$k=1$ to $n$}
		\STATE $R_k^{[i]}=\frac{P_{k-1}^{[i]}+P_k^{[i]}}{2}$
		\STATE $P_k^{[i+1]}=P_{k-1}^{[i]}+\epsilon_k g^k(R_k^{[i]})$
		\ENDFOR
		\STATE $i\to i+1$ 
		\UNTIL{$P_n^{[i]}$ converges.}
		\STATE {Set $j=0$,$T_1^{[0]}=T_2^{[0]}=\cdots=T_n^{[0]}=D_1^{[0]}=D_2^{[0]}=\cdots=D_n^{[0]}=0$}
		\REPEAT
		\STATE $D_{n+1}^{[j]}=f^{\epsilon}(P_n^{[i]})$
		\FOR{$k=n$ to $1$}
		\STATE $T_k^{[j+1]}=f^k(R_k^{[i]})+\frac{\epsilon_k}{2}\partial_{x}g^k(R_k^{[i]})T_k^{[j]}$
		\STATE $B_k^{[j]}=D_{k+1}^{[j]}-T_k^{[j+1]}$
		\STATE $D_k^{[j+1]}=B_k^{[j]}-\frac{\epsilon_k}{2}\partial_{x}g^k(R_k^{[i]})(B_k^{[j]}+D_k^{[j]})$
		\ENDFOR
		\STATE $j\to j+1$ 
		\UNTIL{$D_1^{[j]}$ converges.}		
	\end{algorithmic}
\end{algorithm}

\section{Convergence analysis}\label{sec:convergence-analysis}
\noindent
In this section, we present the convergence analysis of the proposed method. Since our goal is to develop numerical methods to solve ODE systems with a large range of $\bm{\epsilon}$-values, the following assumption appears as a natural prerequisite.  
\begin{assumption}\label{ass-1}
Notice that our method developed in Section \ref{sec:NumSchemeMultiscaleDS} is a first-order method (w.r.t. $\bm{\epsilon}$). We require $f^\epsilon$ and its fluctuation components $f^k$, $k=1,...,n$ are second-order differentiable and are bounded on some closed set $\mathbb{T}^{n}\times\mathcal{K}$, where $\mathcal{K}\subset \mathbb{R}^d$. {  In addition, we assume that the path of the solution  $x(t)$ is in $\mathcal{K}$}.              
\end{assumption}
\begin{remark}
In many cases, $f^\epsilon$ and $f^k$, $k=1,...,n$ are globally defined, which requires $\mathcal{K}=\mathbb{R}^d$.
\end{remark}
First, we prove that the transformed equation \eqref{eqn:non-stiff-nscale-1D} is non-stiff with respect to $\bm\epsilon$, and thus it can be solved by using conventional numerical methods with relatively large time steps.  
\begin{theorem}\label{thm:nonstiff}
Suppose that Assumption \ref{ass-1} is satisfied and $0<\epsilon_{k}<1$, $k=1,...,n$ are sufficiently small. Let $F(t,x)$ denote the right hand side of the ODE system \eqref{eqn:non-stiff-nscale-1D}. Then, we have the following estimate, 
\begin{equation}\label{eqn:unifbdd-dy}
\big|\partial_{t}F(t,x)\big|\leq C_0,
\end{equation}
where $\big|\cdot\big|$ is a vector norm and $C_0$ does not depend on $\epsilon_{k}$, $k=1,...,n$.
\begin{proof}
	According to the definitions \eqref{defi:phik1D}, we have the results,
    \begin{equation}
	\partial_x\Phi^k=I_d+\epsilon_{k}\partial_xg^k, \quad k=1,...,n.
	\end{equation}
	When $\epsilon_{k}$ are sufficiently small, the inverse of $\partial_x\Phi^k$ exists and can be computed through the Neumann series expansion,
	\begin{align}
	(\partial_x\Phi^k)^{-1}&=I_d+\sum_{m=1}^\infty(-\epsilon_{k}\partial_xg^k)^m. 
	\label{Phikinv}
	\end{align}
	Taking the derivative of Eq.\eqref{Phikinv} on both sides with respect to $t$, we obtain
	\begin{align}
	\partial_{t}(\partial_x\Phi^k)^{-1}&=\sum_{m=1}^\infty\partial_{t}(-\epsilon_{k}\partial_xg^k)^m.
	\end{align}
	Moreover, we have the estimates,
	\begin{align}
	\big|\big|(\partial_x\Phi^k)^{-1}-I_d\big|\big|&\leq C\epsilon_{k},\label{DxPhikinvbound}\\
	\big|\big|\partial_{t}(\partial_x\Phi^k)^{-1}\big|\big|&=\big|\big|\sum_{m=1}^\infty\partial_{t}(-\epsilon_{k}\partial_xg^k)^m\big|\big|\leq C,\label{DtPhikinvbound}
	\end{align}
	where $\big|\big|\cdot\big|\big|$ is a matrix norm. At the same time, we have the condition
	\begin{equation}
	\frac{1}{\epsilon_{k}}\partial_{t_k}\Phi^k=\partial_{t_k}g^k=f^k.
	\label{Dt-gkbound}
	\end{equation}
	Therefore, the right hand side of the ODE system \eqref{eqn:non-stiff-nscale-1D} can be re-written as,
	\begin{align}
	F(t,x)&=\prod_{i=1}^{n}(\partial_x\Phi^k)^{-1}f^\epsilon-\sum_{k=1}^n\prod_{i=1}^{k}(\partial_x\Phi^i)^{-1}\frac{1}{\epsilon_k}\partial_{t_k}\Phi^k,\nonumber\\
	&=\prod_{i=1}^{n}(\partial_x\Phi^k)^{-1}\bar{f}^1+\sum_{k=1}^{n-1}\prod_{i=1}^k(\partial_x\Phi^i)^{-1}\big(\prod_{i=k+1}^{n}(\partial_x\Phi^i)^{-1}-I_d\big) f^k,\nonumber\\
	&\equiv J_0+\sum_{k=1}^{n-1}J_k.
	\end{align}
	Taking derivative of $F(t,x)$ with respect to $t$ and using the product rule, we can easily verify that the terms $\partial_{t}J_0$ and $\partial_{t}J_k$, $k=1,...,n-1$ are all $\mathcal{O}(1)$. Thus, the assertion in \ref{eqn:unifbdd-dy} is proved.
\end{proof}
\end{theorem} 
Theorem \ref{thm:nonstiff} shows that the transformed ODE system \eqref{eqn:non-stiff-nscale-1D} is non-stiff, which is then amenable to a standard numerical treatment. As such, we divide the time interval $[0,1]$ by the nodes $t_m=m\Delta t$, $m=0,...,M$, where $\Delta t=1/M$ is the time step and $M$ is a positive integer. For each $m$, $m=1,...,M$, we seek a numerical solution $\hat{y}(t_m)$ to approximate $\tilde{y}(t_m)$, which is the value of the exact solution of the ODE system \eqref{eqn:non-stiff-nscale-1D} at time $t_m$. 

Here, we use the implicit integral midpoint scheme { (\emph{Im2nd})} to solve the ODE system \eqref{eqn:non-stiff-nscale-1D}. { Between two consecutive computational times $t_m$ and $t_{m+1}$, we integrate the differential equation  \eqref{eqn:non-stiff-nscale-1D}  and obtain,
\begin{equation}
\tilde{y}(t_{m+1})=\tilde{y}(t_{m})+\int_{t_{m}}^{t_{m+1}}F(s,\tilde{y}(s))ds.
\end{equation}
Then, we approximate $\tilde{y}(s)$ by an average value and arrive at,
\begin{equation}\label{eqn:im2nd}
	\hat{y}(t_{m+1})=\hat{y}(t_{m})+\int_{t_{m}}^{t_{m+1}}F(s,\frac{\hat{y}(t_{m+1})+\hat{y}(t_{m})}{2})ds.
\end{equation}}
We remark that the numerical solution $\hat{y}(t_{m+1})$ can be computed by some iteration methods, such as the Newton-Raphson method or fixed point iteration method. In this paper, we choose the scheme \eqref{eqn:im2nd} since it preserves certain intrinsic structures in the solution of the original problem; see Section \ref{sec:numericalresults} for more discussions.  

The uniform boundedness of the first-order derivative of $F(t,x)$ (proved in Theorem \ref{thm:nonstiff}) guarantees that our implicit integral midpoint scheme \eqref{eqn:im2nd} has second-order  accuracy. Furthermore, we do not need to decrease the time step $\Delta t$ when $\epsilon_k$, $k=1,...,n$ are small. We summarize the property of the numerical solution $\hat{y}(t_{m})$ into the following lemma.
\begin{lemma}
Let $\tilde{y}(t)$  be the exact solution of the transformed ODE system \eqref{eqn:non-stiff-nscale-1D}. And let 
$\hat{y}(t_m)$, $m=1,...,M$ be the numerical solutions obtained by the scheme \eqref{eqn:im2nd}. Then, we 
have
\begin{align}\label{est:im2nd}
 \big|\hat{y}(t_m)-\tilde{y}(t_m)\big|=C_1(\Delta t)^2,  \quad m=1,...,M,  
\end{align}
where $C_1$ does not depend on $\epsilon_{k}$, $k=1,...,n$.
\end{lemma}
Finally, we analyze the error between the  approximated solution $\bar{\Phi}^n(\hat{y})$ and the 
exact solution $x(t)$ of the original ODE system \eqref{ODEsystem-multiscale}.
	\begin{theorem}\label{thm:mainerror}
	{  Let $T$ denote the final computational time. Suppose Assumption \ref{ass-1} is satisfied and $0<\epsilon_{k}<1$, $k=1,...,n$ are sufficiently small.} For all $t\leq T$, we have the following error estimate
		\begin{align}\label{eqn:final_est}
		\big|x(t)-\bar{\Phi}_t^n(\hat{y}(t))\big|\leq C_2\big(\max_{i=2,\cdots,n}\frac{\epsilon_{i}}{\epsilon_{i-1}}\big)+C_3(\Delta t)^2,
		\end{align}
		where $C_2$ and $C_3$ are generic constants that do not depend on $\epsilon_k$, $k=1,...,n$ and $\Delta t$.
	\begin{proof}
	For any given computational time $t$, we have 	
	\begin{align}
	\big|x(t)-\bar{\Phi}_t^n(\hat{y}(t))\big| \leq
	\big|\bar{\Phi}_t^n(y(t))-\bar{\Phi}_t^n(\tilde{y}(t))\big| + \big|\bar{\Phi}_t^n(\tilde{y}(t))-\bar{\Phi}_t^n(\hat{y}(t))\big|,
	\label{error-splitting}
	\end{align}
	where $y(t)$ and $\tilde{y}(t)$ are the exact solutions of the ODE systems \eqref{master-eq1} and \eqref{eqn:non-stiff-nscale-1D}, respectively, and $\hat{y}(t)$ is the numerical approximation of $\tilde{y}(t)$. We shall estimate the two terms in \eqref{error-splitting} separately. 
	First we can see that,
	\begin{align}
	\dot{y}-\dot{\tilde{y}}=&\big(\prod_{i=1}^{n}(\partial_x\Phi^k)^{-1}f^\epsilon-\sum_{k=1}^n\prod_{i=1}^{k}(\partial_x\Phi^i)^{-1}\sum_{j=1}^{k}\frac{\epsilon_{k}}{\epsilon_j}\partial_{t_j}g^k\big)\nonumber\\&-\big(\prod_{i=1}^{n}(\partial_x\Phi^k)^{-1}f^\epsilon-\sum_{k=1}^n\prod_{i=1}^{k}(\partial_x\Phi^i)^{-1}\partial_{t_k}g^k\Phi^k\big),\nonumber\\
	=&-\sum_{k=1}^n\prod_{i=1}^{k}(\partial_x\Phi^i)^{-1}\sum_{j=1}^{k-1}\frac{\epsilon_{k}}{\epsilon_j}\partial_{t_j}g^k.
	\end{align}
	Using the conditions that $\partial_{t_k}g^k$ are bounded functions (see Eq.\eqref{Dt-gkbound}) and $\big|\big|(\partial_x\Phi^k)^{-1}-I_d\big|\big|\leq C\epsilon_{k}$ (see Eq.\eqref{DxPhikinvbound}), we have 
	\begin{align}
	|\dot{y}-\dot{\tilde{y}}|\leq C(\max_{i=2,\cdots,n}\frac{\epsilon_{i}}{\epsilon_{i-1}}).
	\end{align}
	Hence for any time $t\leq T$, we have, 
	\begin{align}\label{est:UAapprox}
	|y(t)-\tilde{y}(t)|\leq C_T (\max_{i=2,\cdots,n}\frac{\epsilon_{i}}{\epsilon_{i-1}}),
	\end{align}
	where the constant $C_T=\mathcal{O}(T)$. 
	Applying the chain rule for $\bar{\Phi}^n$, we obtain 
	\begin{align}
	\partial_x\bar{\Phi}^n=\partial_x\big(\Phi^n\circ\Phi^{n-1}\circ\cdots\circ\Phi^1\big) =\prod_{i=1}^{n}(I_d+\epsilon_{i}\partial_x g^i).
	\label{DxPhinBounded}
	\end{align}
	So at any time $t$, $\bar{\Phi}^n$ is Lipschitz in the variable $x$ and the Lipschitz constant is uniformly bounded when $\epsilon_i$ are small enough. 
	Finally, combining the estimates \eqref{est:im2nd}, \eqref{est:UAapprox} and \eqref{DxPhinBounded}, we prove the statement in Theorem \ref{thm:mainerror}.
    \end{proof}
\end{theorem}

When the multiscale parameters are well-separated; see \eqref{WellSeparated-multiscale}, the first term in the error estimate \eqref{eqn:final_est} is negligible, thus our scheme has a second-order accuracy with respect to $\Delta t$. Although the above convergence analysis relies on the scale separation assumption, we can relax this assumption in some special cases. For example, when two scales collapse, i.e. $\epsilon_{n} = c \epsilon_{n-1}$, we can treat these two scales as a single scale and we can modify the mean function defined in \eqref{mean-fn} accordingly. More specifically, if $c =\frac{m_1}{m_2}$ is a rational number, and $h(y_1, y_2)$ is a doubly periodic function in $y_1$ and $y_2$ with period $[1,1]$, then
$h(\frac{t}{\epsilon_{n-1}}, \frac{t}{\epsilon_{n}}) = h(\frac{t}{\epsilon_{n-1}}, 
\frac{m_2 t} {m_1 \epsilon_{n-1}})$ is a periodic function of $t$ with a period $m_1\epsilon_{n-1}$.  Thus, the mean function defined in \eqref{mean-fn} can be modified accordingly as follows
\begin{equation}
\frac{1}{m_1} \int_0^{m_1} f(t_1, t_2, ..., t_{n-2}, s, \frac{m_2}{m_1} s, x) ds.
\end{equation}
When $c$ is an irrational number, it is easy to show that the time average of $h(\frac{t}{\epsilon_{n-1}}, \frac{t}{\epsilon_{n}})$ will converge to the area in $(y_1, y_2)\in \mathbb{T}^{2}$ and we can modify the definition of the mean function as follows
\begin{equation}
\int_{[0,1]^2} f(t_1, t_2, ..., t_{n-2}, s_1, s_2 , x) ds_1 ds_2 .
\end{equation}
With the above modification of the mean function, we can still prove the main results stated in this section. Our numerical results to be presented later also confirm that our method works equally well in the case when two scales collapse. 

More general case can be considered as well if we have $m$ number of collapsed scales, i.e.
$\epsilon_k = c_1 \epsilon_{k+1} = c_2 \epsilon_{k+2} = ...= c_{m-1} \epsilon_{k+m-1}$.
In this case, we should consider these m-scales simultaneously and modify the definition of the mean function by using the time averaging technique for multiple scales discussed in \cite{hou1988homogenization}. We will not present the more general case in this paper and will leave it to our future work.

\section{Numerical results} \label{sec:numericalresults}
\noindent 
In this section, we present several numerical experiments to illustrate the efficiency of our method and confirm the convergence analysis. { The \emph{Im2nd} solution, or direct \emph{Im2nd} solution, refers to the numerical solution obtained by solving the original problem \eqref{ODEsystem-multiscale} using the scheme \eqref{eqn:im2nd}}. Moreover, we use the \emph{UA} solution to denote the numerical solution obtained by our method. 
In our method, we use the scheme \eqref{eqn:im2nd} to solve Eq.\eqref{eqn:non-stiff-nscale-1D} and 
the Algorithm \ref{alg:nscale} to compute necessary quantities in the Eq.\eqref{eqn:non-stiff-nscale-1D}.

\subsection{An ODE system with three separated time scales} 
\noindent 
The H\'{e}non-Heiles system \cite{henon1964applicability} is undoubtedly one of the most paradigmatic model potentials for time-independent Hamiltonian systems with two degrees of freedom, which is frequently used to 
describe the motion of stars around a galactic center. We consider a generalization of the original H\'{e}non-Heiles system in three degrees of freedom. The relative equilibria and bifurcations of this model were 
studied in \cite{ferrer1998henon}. 

Here, we assume the Hamiltonian of the three dimensional H\'{e}non-Heiles system is parameterized by $\epsilon_1$ and $\epsilon_2$ and has the following form  

\begin{equation}
H(\textbf{p},\textbf{q})=\frac{p_1^2}{2\epsilon_2}+\frac{q_1^2}{2\epsilon_2}+\frac{p_2^2}{2\epsilon_1}+\frac{q_2^2}{2\epsilon_1}+\frac{p_3^2}{2}+\frac{q_3^2}{2}+q_1^2q_2-\frac{1}{3}q_2^3+q_2^2q_3-\frac{1}{3}q_3^3.
\label{Hamiltonian-3scales}
\end{equation}
where $\textbf{p}=(p_1,p_2,p_3)^{T}$ and $\textbf{q}=(q_1,q_2,q_3)^{T}$. One can obtain the evolution equation for the Hamiltonian in \eqref{Hamiltonian-3scales} as follows, 
\begin{equation}
\frac{d \textbf{p}}{dt}=-\frac{\partial H}{\partial \textbf{q}},  \quad \quad \quad 
\frac{d \textbf{q}}{dt}= \frac{\partial H}{\partial \textbf{p}}.
\label{Hamiltonian-3scales-ODE}
\end{equation}  
When $0<\epsilon_2\ll\epsilon_1\ll1$, the ODE system \eqref{Hamiltonian-3scales-ODE} becomes a three-scale problem and the solutions (e.g., $p_1$ and $q_1$) are highly oscillatory. Let us carry out a change of variables,  

\begin{equation}
\begin{cases}
w_1=\cos(\frac{t}{\epsilon_2})q_1-\sin(\frac{t}{\epsilon_2})p_1,\\
w_2=\sin(\frac{t}{\epsilon_2})q_1+\cos(\frac{t}{\epsilon_2})p_1,\\
w_3=\cos(\frac{t}{\epsilon_1})q_2-\sin(\frac{t}{\epsilon_1})p_2,\\
w_4=\sin(\frac{t}{\epsilon_1})q_2+\cos(\frac{t}{\epsilon_1})p_2,\\
w_5=q_3,\\
w_6=p_3.\label{changevariable-3scales}
\end{cases}
\end{equation}
Then, $(w_1,...,w_6)^{T}$ satisfies the following ODE system
\begin{equation}\label{ODE-3scales}
\begin{small}
\begin{cases}
\dot{w}_1=2 \sin t_2(w_1\cos t_2+w_2\sin t_2)(w_3\cos t_1+w_4\sin t_1),\\
\dot{w}_2=-2\cos t_2(w_1\cos t_2+w_2\sin t_2)(w_3\cos t_1+w_4\sin t_1),\\
\dot{w}_3=\sin t_1\big(2(w_3\cos t_1+w_4\sin t_1)w_5+(w_1\cos t_2+w_2\sin t_2)^2-(w_3\cos t_1+w_4\sin t_1)^2\big),\\
\dot{w}_4=-\cos t_1\big(2(w_3\cos t_1+w_4\sin t_1)w_5+(w_1\cos t_2+w_2\sin t_2)^2-(w_3\cos t_1+w_4\sin t_1)^2\big),\\
\dot{w}_5=w_6,\\
\dot{w}_6=w_5^2-w_5-(w_3\cos t_1+w_4\sin t_1)^2,
\end{cases}
\end{small}
\end{equation} 
where $ t_2=\frac{t}{\epsilon_2}$ and $ t_1=\frac{t}{\epsilon_1}$. 
Notice that the right-hand side of the ODE system \eqref{ODE-3scales} involves trigonometric functions and simple polynomials. Thus, all the integrals in our numerical schemes can be pre-computed analytically. We have implemented these computations with the software Mathematica.  
 
\textbf{Verification of the convergence analysis.} We compare the error between the numerical solution obtained by our method and the reference solution.  The initial value is  $(w_1(0),...,w_6(0))^{T}=(0.12,0.12,0.12,0.12,0.12,0.12)^{T}$.  The reference solution is obtained by Matlab \texttt{ode45} function applied to the ODE system \eqref{ODE-3scales}, where the time step is $\Delta t=0.0001$. To implement our method, we choose the second-order implicit integral midpoint scheme to integrate the non-stiff problem \eqref{eqn:non-stiff-nscale-1D}. We find that the iteration loops needed in Algorithm \ref{alg:nscale} is about $3-10$ times, where the convergence threshold is set to be $10^{-14}$.  

\begin{figure}[htbp]
	\centering
	\subfigure[]{
		\includegraphics[width=0.49\linewidth]{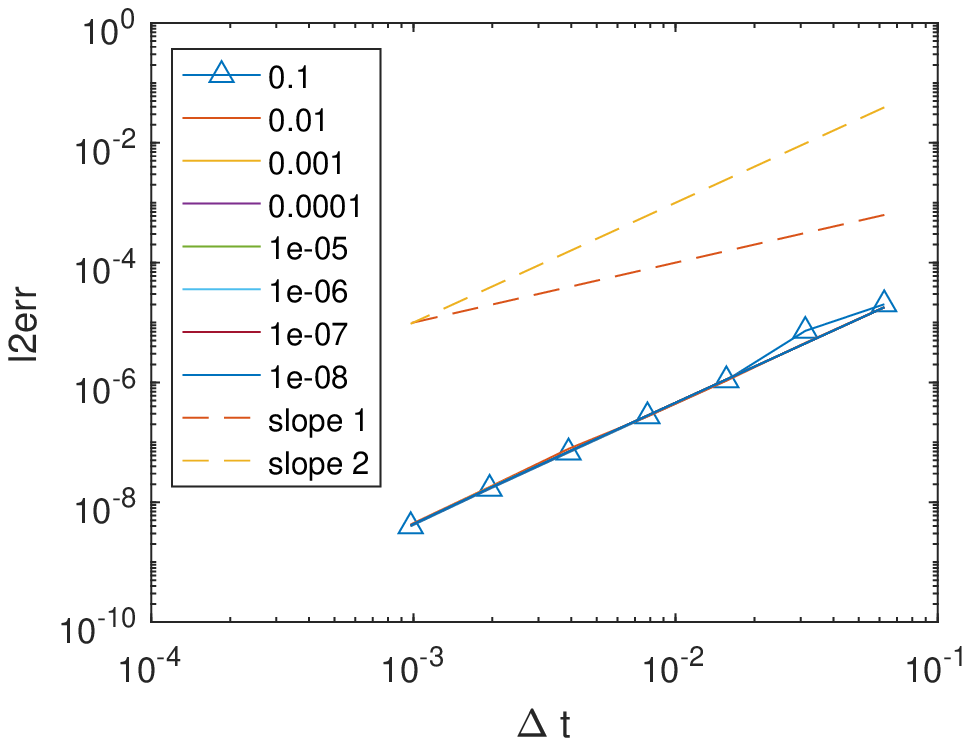}
		\label{fig:HHmodel_errT1}
	}%
	\subfigure[]{
		\includegraphics[width=0.49\linewidth]{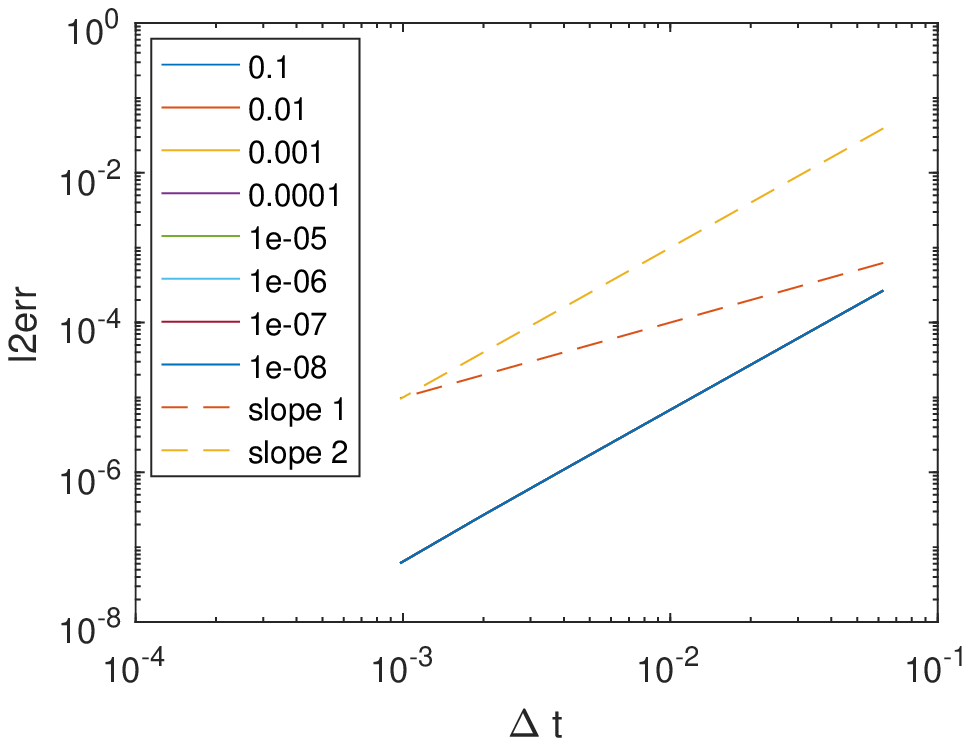}
		\label{fig:HHmodel_errv3}
	} 
	\caption{Error as a function of $\Delta t$ for $\epsilon_1=10^{-k}$, $k=1,...,8$. Left: $\epsilon_2=\epsilon_1^2$. Right:  $\epsilon_2=0.8\epsilon_1$.}
    \label{fig:3scale_errordt1}
	\end{figure}

In Fig.\ref{fig:HHmodel_errT1} and Fig.\ref{fig:HHmodel_errv3}, we show the error as a function of $\Delta t$ for different values of $\epsilon_1$ with $\epsilon_2=\epsilon_1^2$ and $\epsilon_2=0.8\epsilon_1$, respectively. 
The magnitude of the Hamiltonian \eqref{Hamiltonian-3scales} is about $\frac{0.0144}{\epsilon_1^2}$. We observe a second-order convergence rate with respect to $\Delta t$ in our method. Most importantly, the error is independent of $\epsilon_1$ and $\epsilon_2$ as shown in Fig.\ref{fig:HHmodel_errT1} and Fig.\ref{fig:HHmodel_errv3}, where the curves for different values of $\epsilon_1$ and $\epsilon_2$ are nearly identical. This confirms that our scheme is uniform
accurate, which does not depend on $\epsilon_1$ and $\epsilon_2$.  Notice that the numerical solution obtained by $\Delta t=0.1$ is accurate enough to maintain an error of the order $10^{-4}$ and error of the Hamiltonian at the order of $10^{-3}$, uniformly in $\epsilon_1$ and $\epsilon_2$. For the Euler method, it is impossible since the 
ODE system associated with the Hamiltonian \eqref{Hamiltonian-3scales} becomes severely stiff when $\epsilon_1$ and $\epsilon_2$ become small. The numerical results for the Euler method were not shown here. 

In the derivation of the numerical method and the convergence analysis, we assume the time-scales are well-separated, i.e., $0<\epsilon_n\ll\epsilon_{n-1}\ll\cdot\cdot\cdots\ll\epsilon_1\ll1$.  In Fig.\ref{fig:HHmodel_errv3}, we show the error as a function of $\Delta t$ for difference values of $\epsilon_1$ and  $\epsilon_2=0.8\epsilon_1$ at $T=1$. It is shown  that our scheme still has a uniform accuracy, which does not depend on separation between $\epsilon_1$ and $\epsilon_2$.

\begin{figure}[htbp]
	\centering
	\subfigure[]{
		\includegraphics[width=0.49\linewidth]{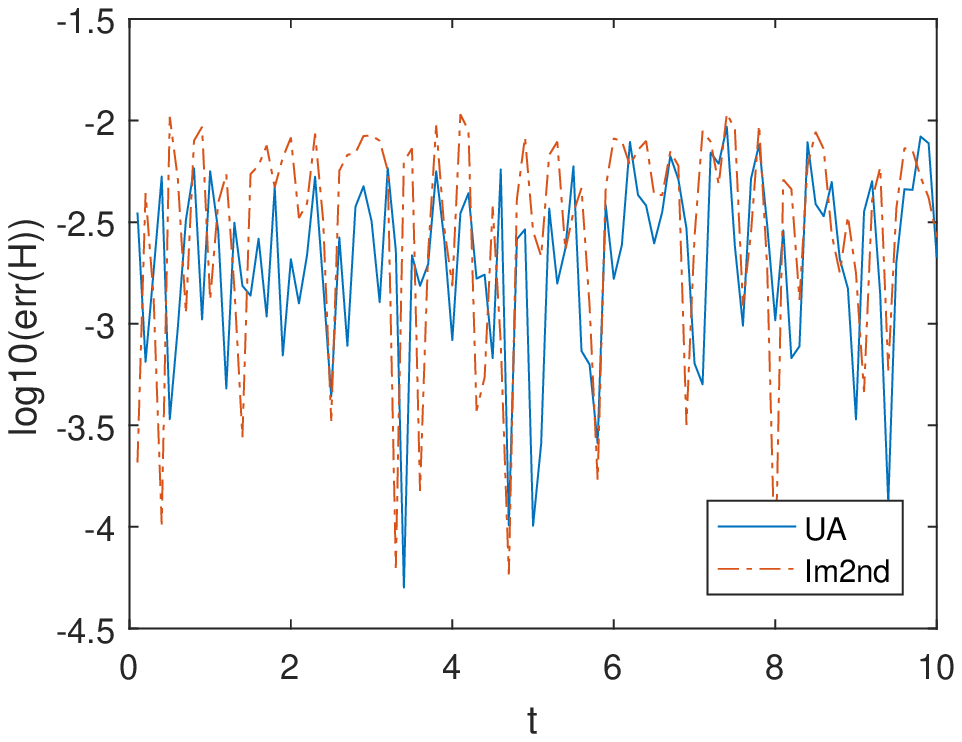}
		\label{fig:eps1e-1}
	}%
	\subfigure[]{
		\includegraphics[width=0.49\linewidth]{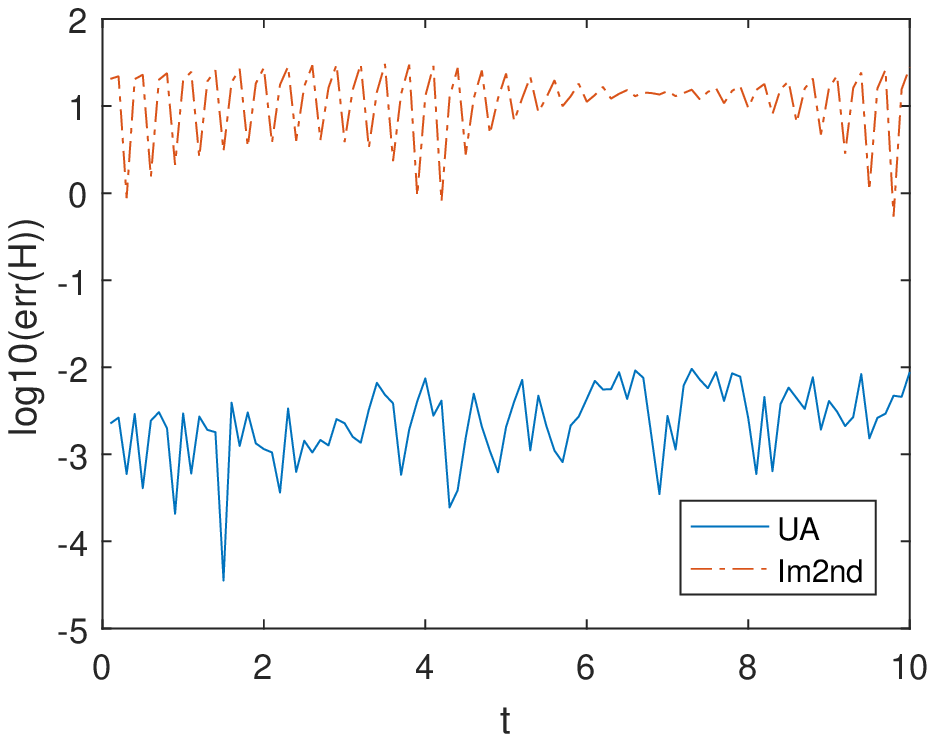}
		\label{fig:eps1e-3}
	}%
	\caption{Evolution of the error of the Hamiltonian in the three-scale ODE system with $\epsilon_{2}=\epsilon_{1}^2$. Left: $\epsilon_1=0.1$. Right: $\epsilon_1=0.001$. $\Delta t=0.1$.}
	\label{fig:H_err_3scale}
\end{figure}

Let us now verify that our method preserves the Hamiltonian of the system.  In Fig.\ref{fig:H_err_3scale}, we plot the evolution of the error of the Hamiltonian \eqref{Hamiltonian-3scales} for different $\epsilon_1$ and $\epsilon_2$, where $\epsilon_{2}=\epsilon_{1}^2$. We find that when $\epsilon_{1}$ is relatively large, e.g., $\epsilon_{1}=0.1$, the implicit integral midpoint scheme \eqref{eqn:im2nd} with 
a large time step $\Delta t=0.1$ and our method give similar and accurate results in computing the Hamiltonian; see Fig.\ref{fig:eps1e-1}. However, when $\epsilon_{1}$ is small, e.g., $\epsilon_{1}=0.001$, the ODE system associated with the Hamiltonian \eqref{Hamiltonian-3scales} becomes very stiff. 
Directly using the scheme \eqref{eqn:im2nd} will lose accuracy if the time step is not small enough. While our method with a large time step still maintains the same accuracy; see Fig.\ref{fig:eps1e-3}. This result again confirms that our scheme has a uniform accuracy in computing an ODE system with multiple time-scales, especially when the systems are stiff. 

\textbf{Verification of the non-stiffness of the transformed equation \eqref{eqn:non-stiff-nscale-1D}}.
In Theorem \ref{thm:nonstiff}, we proved that under certain assumptions the transformed ODE system \eqref{eqn:non-stiff-nscale-1D} is non-stiff. Here, we shall verify this statement numerically by 
solving the three-scale ODE system \eqref{ODE-3scales} with an initial value $(w_1(0),...,w_6(0))^{T}=(0.12,0.12,0.12,0.12,0.12,0.12)^{T}$. We compute the maps 
$\Phi^1$ and $\Phi^2$ in our method with different $t_1$ and $t_2$ for $(y_1(0),...,y_6(0))^{T}=(0.20,0.20,0.20,0.20,0.20,0.20)^{T}$. In addition, we record the quantities 
{ $P_2 = \Phi_{ t_2}^2(\Phi_{ t_1}^1(y))$}, $T_1=\frac{1}{\epsilon_1}\partial_{ t_1}\Phi^1_{ t_1}(y)$, 
$T_2=\frac{1}{\epsilon_2}\partial_{ t_2}\Phi^2_{ t_2}(\Phi^1_ {t_1}(y))$, and $D_1= \dot{y}$. Recall that these quantities were defined in Section \ref{sec:Implentmentation}, especially 
in the Algorithm \ref{alg:nscale}.

\begin{figure}[htbp]
	\centering
	\includegraphics[width=1.0\linewidth]{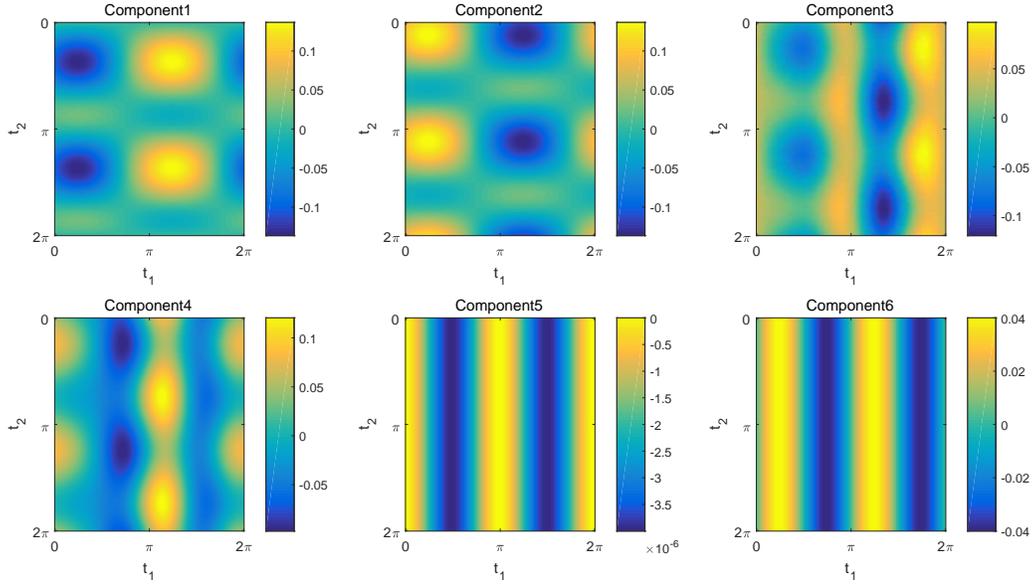}
	\caption{Six components of { $f(t_1,t_2,y)-D_1$}. Here $\epsilon_1=0.0001$ and $\epsilon_2=\epsilon_1^2$.}
	\label{fig:HH_1e-4_1e-8_f-D1}
\end{figure}

In Fig.\ref{fig:HH_1e-4_1e-8_f-D1}, we show the six components of the quantities $f(y)-D_1$ as functions of $ t_1$ and $ t_2$ (i.e., the six components of right hand side of the ODE system \eqref{ODE-3scales} minus their numerical counterparts), where $\epsilon_1=0.0001$ and $\epsilon_2=\epsilon_1^2$. 
In this example, $\max|f(y)|$ is $\mathcal{O}(1)$. We also compute the cases when $\epsilon_1=0.01$ and $\epsilon_1=0.001$ with $\epsilon_2=\epsilon_1^2$ and find the patterns of the six components of the quantities $\dot{y}-D_1$ remain almost the same as the Fig.\ref{fig:HH_1e-4_1e-8_f-D1}. Thus, we do not show them here. 

\begin{figure}[htbp]
	\centering
	\includegraphics[width=1.0\linewidth]{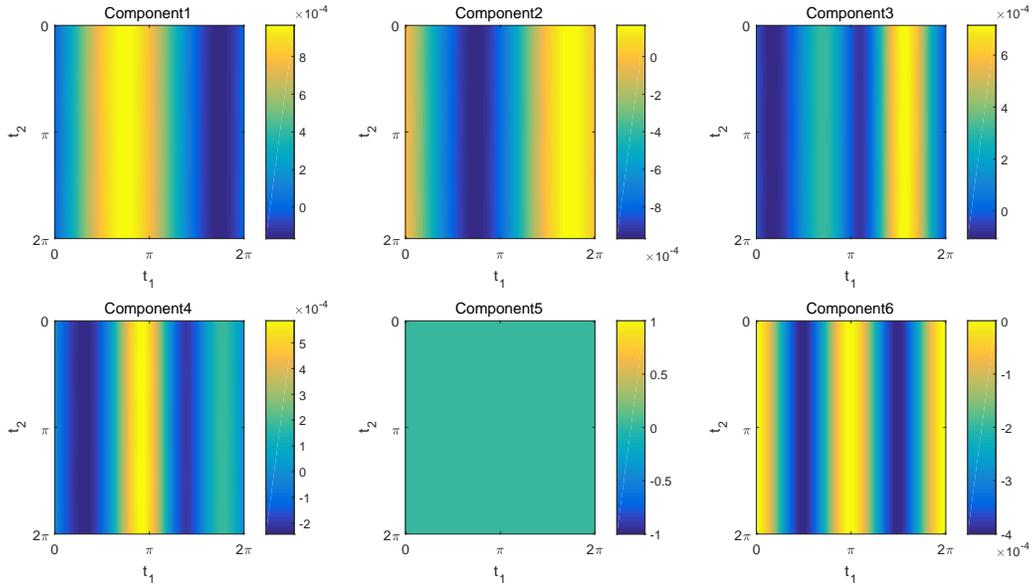}
	\caption{$\Phi_{ t_2}^2(\Phi_{ t_1}^1(y))-y$ when $\epsilon_1=0.01$ and $\epsilon_2=\epsilon_1^2$.}
	\label{fig:HH_1e-2_1e-4_P2}
\end{figure}

\begin{figure}[htbp]
	\centering
	\includegraphics[width=1.0\linewidth]{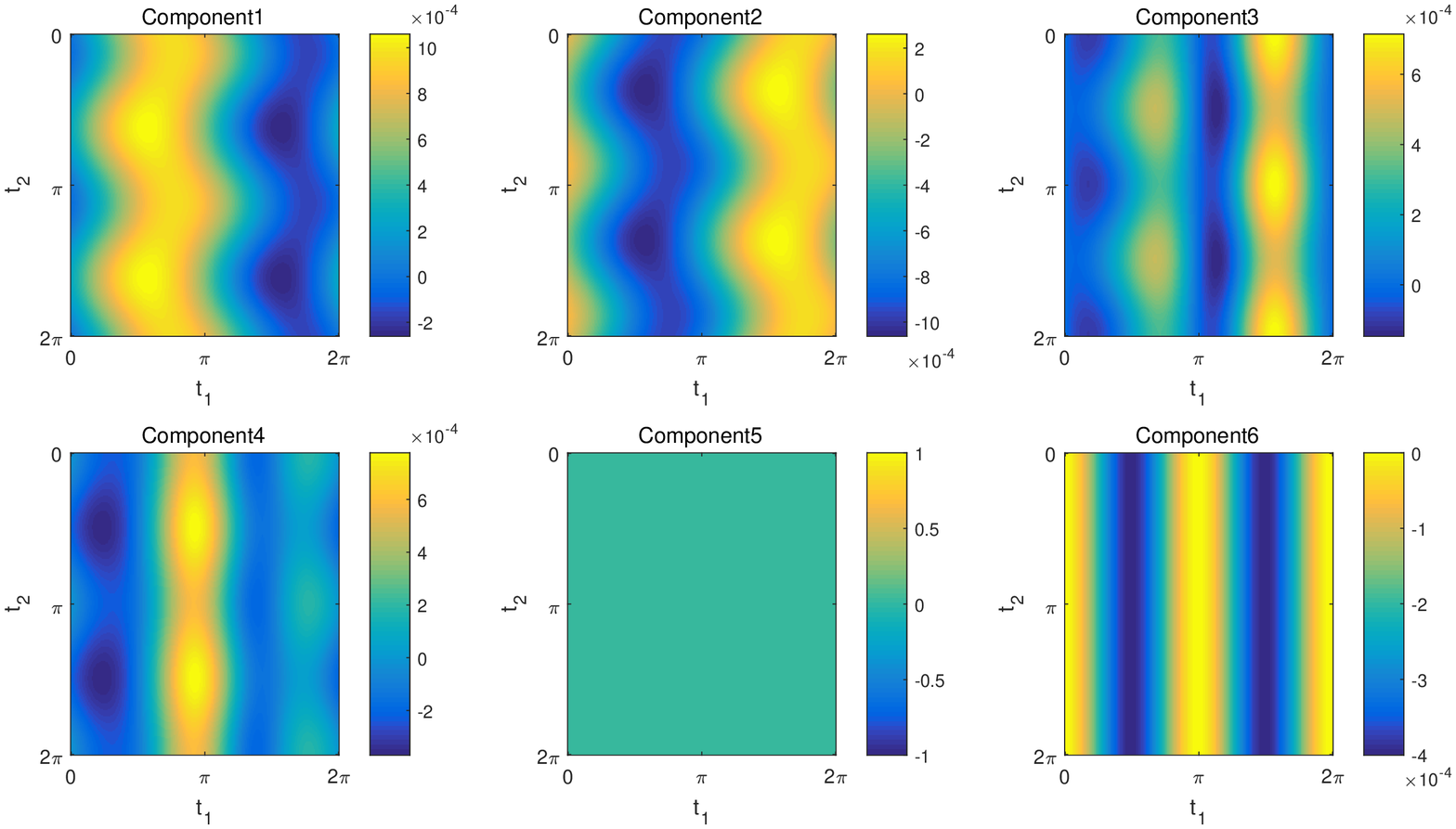}
	\caption{$\Phi_{ t_2}^2(\Phi_{ t_1}^1(y))-y$ when $\epsilon_1=0.01$ and $\epsilon_2=\frac{\epsilon_1}{2}$.}
	\label{fig:HH_1e-2_5e-3_P2}
\end{figure}

In Fig.\ref{fig:HH_1e-2_1e-4_P2}, we show the magnitude of the quantity $\Phi_{ t_2}^2(\Phi_{ t_1}^1(y))-y$ as functions of $ t_1$ and $ t_2$ when $\epsilon_1=0.01$ and $\epsilon_2=\epsilon_1^2$. The results for $\Phi_{ t_2}^2(\Phi_{ t_1}^1(y))-y$ when $\epsilon_1=0.01$ and  $\epsilon_2=\frac{\epsilon_1}{2}$ are shown in Fig.\ref{fig:HH_1e-2_5e-3_P2}. One can see that when there is no scale separation $\Phi_{ t_2}^2(\Phi_{ t_1}^1(y))$ are fluctuating along $ t_2$ direction. This result provides numerical confirmation of our derivation. We can write $\Phi_{ t_2}^2(\Phi_{ t_1}^1(y))$ explicitly out as, 
\begin{align}
\Phi_{ t_2}^2(\Phi_{ t_1}^1(y)) =\Phi_{ t_1}^1(y)+\epsilon_2g^2( t_2, t_1,\Phi_{ t_1}^1(y)) 
 =y+\epsilon_{1}g^1( t_1,y)+\epsilon_2g^2( t_2, t_1,y+\epsilon_{1}g^1( t_1,y)).
\end{align}
Then given $y$, when $\epsilon_1$ and $\epsilon_2$ are close, $\epsilon_2g^2( t_2, t_1,y+\epsilon_{1}g^1( t_1,y))$ is comparable to $\epsilon_{1}g^1( t_1,y)$. 

 
\begin{figure}[htbp]
	\centering
	\includegraphics[width=1.0\linewidth]{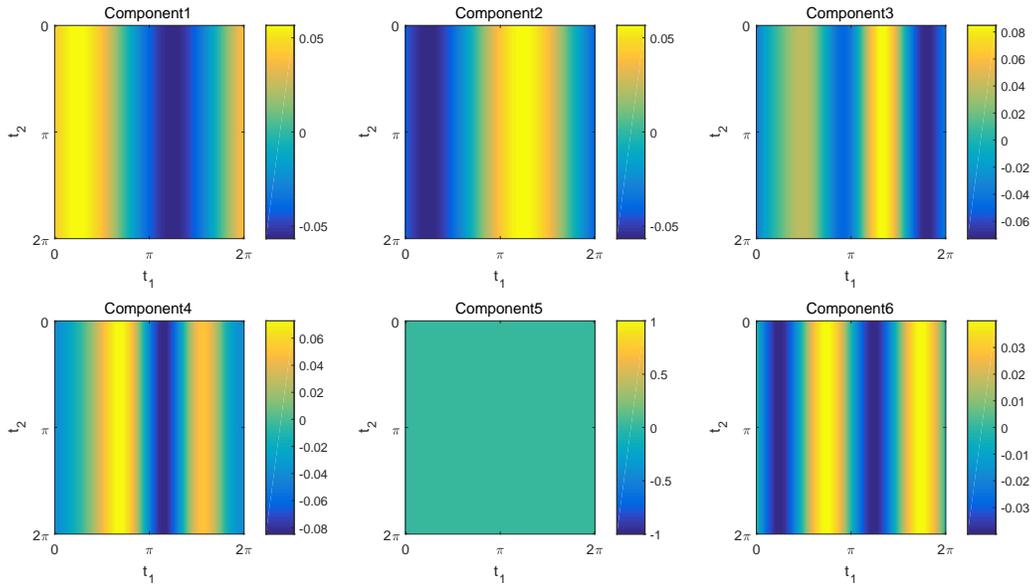}
	\caption{$T_1$ when $\epsilon_1=0.0001$ and $\epsilon_2=\epsilon_1^2$.}
	\label{fig:HH_1e-4_1e-8_T1}
\end{figure}

\begin{figure}[htbp]
	\centering
	\includegraphics[width=1.0\linewidth]{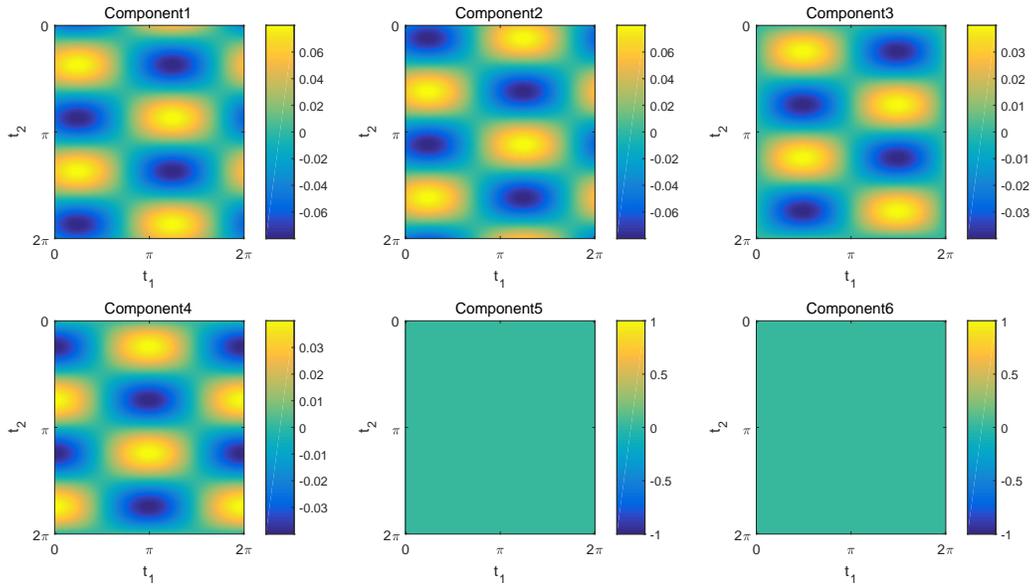}
	\caption{$T_2$ when $\epsilon_1=0.0001$ and $\epsilon_2=\epsilon_1^2$.}
	\label{fig:HH_1e-4_1e-8_T2}
\end{figure}

In Fig.\ref{fig:HH_1e-4_1e-8_T1} and Fig.\ref{fig:HH_1e-4_1e-8_T2}, we show the magnitude of the quantities $T_1$ and $T_2$ as functions of $ t_1$ and $ t_2$, when $\epsilon_1=0.0001$ and $\epsilon_2=\epsilon_1^2$, respectively.  We also compute the  quantities $T_1$ and $T_2$ when $\epsilon_1$ changes (e.g. $\epsilon_1=0.01$ and $\epsilon_1=0.001$) with $\epsilon_2=\epsilon_1^2$ and find the time derivatives of $T_1$ and $T_2$ remain almost the same magnitude as that shown in Fig.\ref{fig:HH_1e-4_1e-8_T1} and Fig.\ref{fig:HH_1e-4_1e-8_T2}. This implies the implicit iterative scheme Eq.\eqref{eqn:iterative-nscale} for Eq.\eqref{midpoint-Phik} keeps the magnitude of the non-midpoint setting, $\partial_{t}\Phi^k=\sum_{i=1}^{k}\frac{\epsilon_{k}}{\epsilon_{i}}\partial_{t_i}g^k$. Notice that component 5 and 6 do not depend on $ t_2$. This is due to the fact that $f_5$ and $f_6$ are independent of $ t_2$. $T_1$ is independent of $ t_2$ for difference choice of  $\epsilon_1$ and $\epsilon_2$, this is due to the definition of $T_1$; see Fig.\ref{fig:HH_1e-4_1e-8_T1}.

\subsection{An ODE system with four separated time scales}
\noindent  
To further study the performance of our method, we mimic the formulation of the three dimensional 
H\'{e}non-Heiles system and generate a Hamiltonian system with four time scales. The Hamiltonian is given by,
\begin{align}
H(\textbf{q},\textbf{p})=&\frac{p_1^2}{2\epsilon_3}+\frac{q_1^2}{2\epsilon_3}+\frac{p_2^2}{2\epsilon_2}+\frac{q_2^2}{2\epsilon_2}+\frac{p_3^2}{2\epsilon_1}+\frac{q_3^2}{2\epsilon_1}+\frac{p_4^2}{2}+\frac{q_4^2}{2}\nonumber\\&+q_1^2q_2+q_2^2q_3+q_3^2q_4-\frac{1}{3}q_2^3-\frac{1}{3}q_3^3-\frac{1}{3}q_4^3,
\label{Hamiltonian-4scales}
\end{align}
where $\textbf{p}=(p_1,p_2,p_3,p_4)^{T}$ and $\textbf{q}=(q_1,q_2,q_3,q_4)^{T}$. 
The Hamiltonian in \eqref{Hamiltonian-4scales} is parameterized by $\epsilon_1$, $\epsilon_2$, and $\epsilon_3$. One can obtain the evolution equation for the Hamiltonian in \eqref{Hamiltonian-4scales} 
by using the relation $\frac{d \textbf{p}}{dt}=-\frac{\partial H}{\partial \textbf{q}}$ and $\frac{d \textbf{q}}{dt}= \frac{\partial H}{\partial \textbf{p}}$. When $0<\epsilon_3\ll\epsilon_2\ll\epsilon_1\ll1$, the associated evolution equation of the Hamiltonian in \eqref{Hamiltonian-4scales} becomes a four-scale ODE system. 
Since the formulation of the change of variables and derivation of the transformed ODE system are  standard (similar as we did in Eqns.\eqref{changevariable-3scales} and \eqref{ODE-3scales}), we do not show them here.  
 
\begin{figure}[htbp]
	\centering
	\subfigure[]{
		\includegraphics[width=0.49\linewidth]{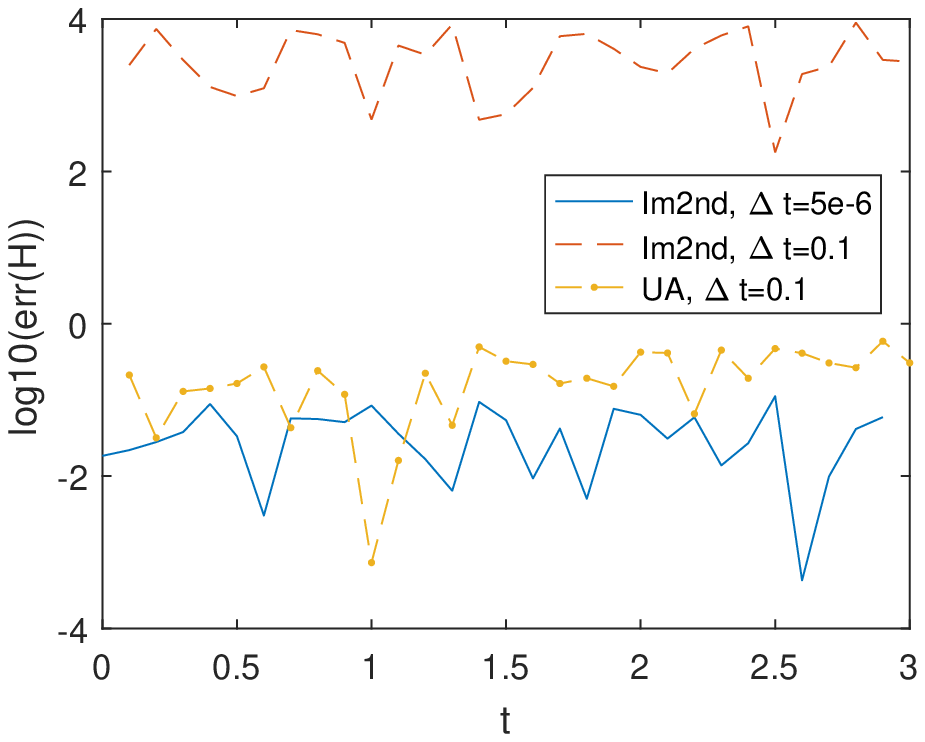}		
		\label{fig:hh4hamiltonian1e-6}
	}%
	\subfigure[]{
		\includegraphics[width=0.49\linewidth]{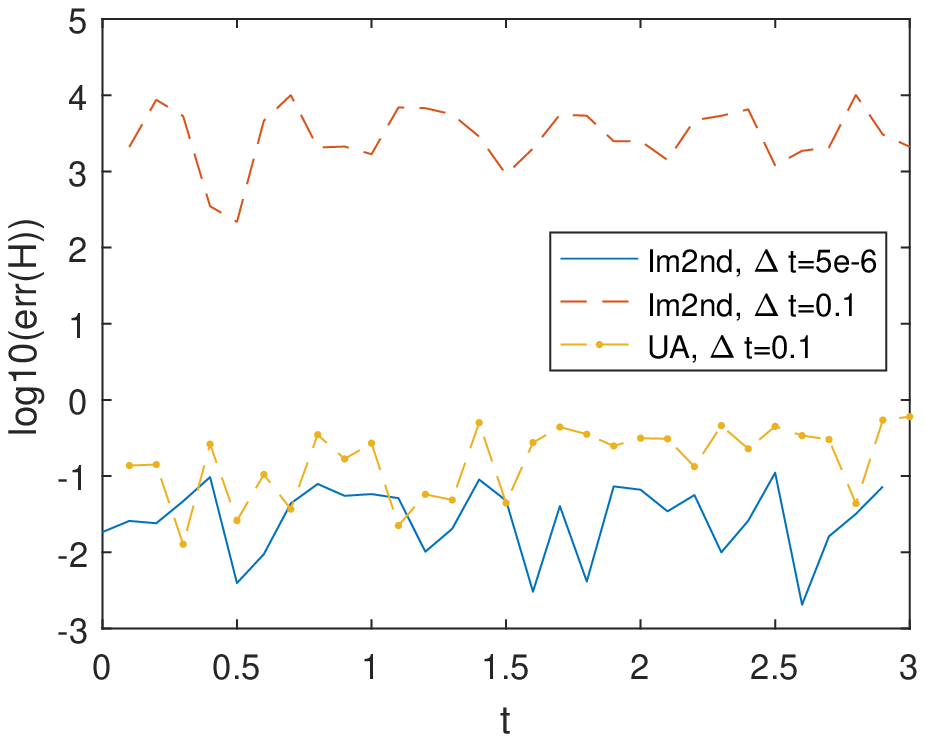}
		\label{fig:hh4hamiltonian1e-6v2}
	}
    \caption{Evolution of the error of the Hamiltonian in the four-scale ODE system. 
	Left:  $(\epsilon_1,\epsilon_2,\epsilon_3)=(10^{-3},11\times 10^{-5}, 3\times 10^{-6})$. 
	Right: $(\epsilon_1,\epsilon_2,\epsilon_3)=(0.5,11\times 10^{-6}, 3\times 10^{-6})$.}
    \label{fig:hh4hamiltonian}
\end{figure}

\textbf{Verification of the convergence analysis.}   Let us now check that our method preserves the Hamiltonian of the system.  In Fig.\ref{fig:hh4hamiltonian}, we present the error of the Hamiltonian \eqref{Hamiltonian-4scales} as a function of time for different $\epsilon_i$, $i=1,2,3$. 
In Fig.\ref{fig:hh4hamiltonian1e-6}, we set the multiscale time-scale parameters to be $(\epsilon_1,\epsilon_2,\epsilon_3)=(10^{-3},11\times 10^{-5}, 3\times 10^{-6})$ and the initial values are $w_i(0)=0.44$ for $i=1,\cdots,8$. In Fig.\ref{fig:hh4hamiltonian1e-6v2}, we set $(\epsilon_1,\epsilon_2,\epsilon_3)=(0.5,11\times 10^{-6}, 3\times 10^{-6})$ and the initial values are the same. Notice that the smallest period of the ODE system is about $2\pi\times\epsilon_3\approx 1.885\times10^{-5}$. To resolve the oscillation in the solution, we choose the time step for the fine-scale ODE solver to be $\Delta t=5\times 10^{-6}$. To implement our method, we choose the time step to be $\Delta t=  10^{-1}$. 

Numerical results in Fig.\ref{fig:hh4hamiltonian}~show that: (1) directly using the implicit integral midpoint scheme \eqref{eqn:im2nd} with a coarse time step $\Delta t= 10^{-1}$ gives wrong results; (2) our method with the same coarse time step $\Delta t= 10^{-1}$ gives an accurate result that is comparable to that using the scheme \eqref{eqn:im2nd} with a very fine time step  $\Delta t= 10^{-6}$. This comparison again confirms that our scheme has a uniform accuracy in computing ODE system with multiple time-scales. In addition, from the results in Fig.\ref{fig:hh4hamiltonian1e-6v2}, where $(\epsilon_1,\epsilon_2,\epsilon_3)=(0.5,11\times 10^{-6}, 3\times 10^{-6})$, we can see that when the time scales are not well separated ( i.e., $\epsilon_3$ is close to $\epsilon_2$ and $\epsilon_{1}$ is close to $1$), our numerical method still gives an excellent performance. 

In terms of the computational time, our method takes $0.163$ and $0.367$ seconds to compute the result shown in Fig.\ref{fig:hh4hamiltonian1e-6} and Fig.\ref{fig:hh4hamiltonian1e-6v2}, respectively. While the direct \emph{Im2nd} method with $\Delta t= 5\times10^{-6}$ takes  about $10.48$ seconds for both. 
Thus, our method achieves a $20\sim60X$ speedup over the conventional ODE solver in this example. 
Moreover, our method provides uniform accurate results for different values of $\epsilon_i$, $i=1,2,3$. 
The conventional ODE solvers, such as the direct \emph{Im2nd} method, require to choose finer time steps when we decrease 
$\epsilon_i$, $i=1,2,3$.  Therefore, it is expected that a higher speedup will be achieved when we need to solve
a multiscale ODE system with much smaller  $\epsilon_i$, $i=1,2,3$. 
 
\begin{figure}[htbp]
	\centering
	\includegraphics[width=0.6\linewidth]{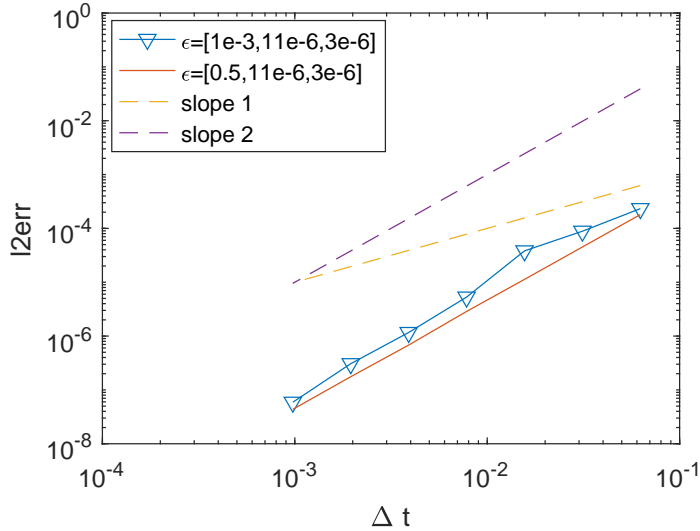}
	\caption{Error as a function of $\Delta t$ for two sets of multiscale parameters at time $T=3$.}
	\label{fig:HH4model_err}
\end{figure}

We also investigate the convergence rate of our method with respect to the time step. 
In Fig.\ref{fig:HH4model_err}, we show the error as a function of $\Delta t$ for two sets of values of $\epsilon_1$, $\epsilon_2$, and $\epsilon_3$  at time $T=3$, respectively. The intial values $x_i(0)$, $i=1,...,8$ and values of $(\epsilon_1,\epsilon_2,\epsilon_3)$ are the same as before. The reference solution is obtained by Matlab \texttt{ode45} function applied to the ODE system, where the time step is $\Delta t=5\times10^{-7}$. { We observe a second-order convergence rate with respect to $\Delta t$ in our method, which verifies the error estimate in Theorem \ref{thm:mainerror}.} Most importantly, the error is independent of $\epsilon_1$, $\epsilon_2$, and $\epsilon_3$ as shown in Fig.\ref{fig:HH4model_err}. 

\begin{figure}[htbp]
	\centering
	\subfigure[]{
		\includegraphics[width=0.49\linewidth]{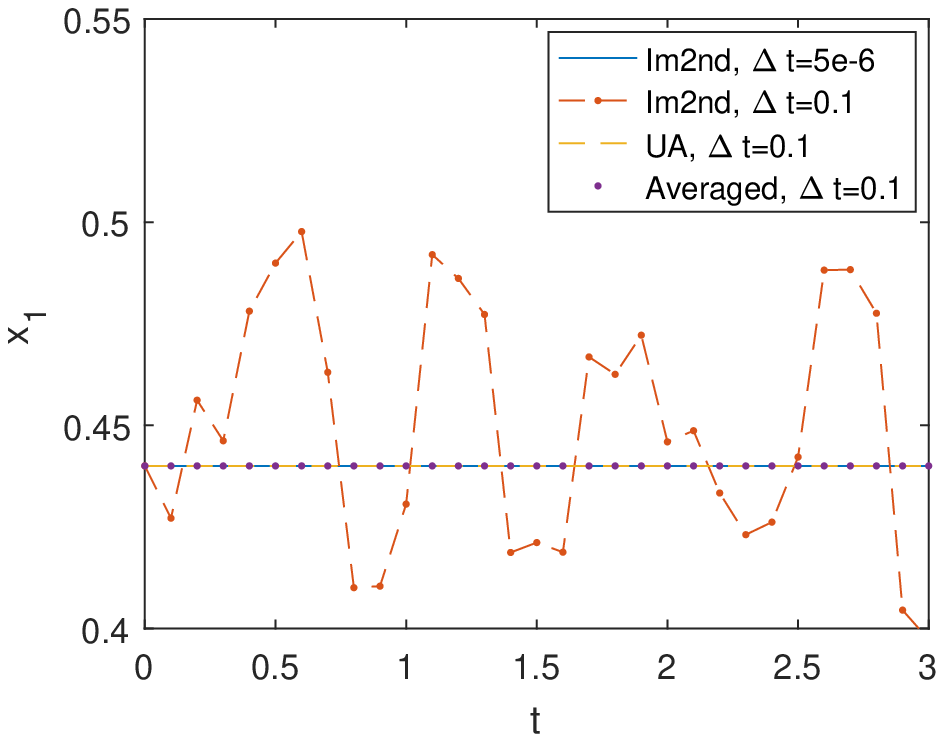}
		\label{fig:x1_avg}
	}%
	\subfigure[]{
		\includegraphics[width=0.49\linewidth]{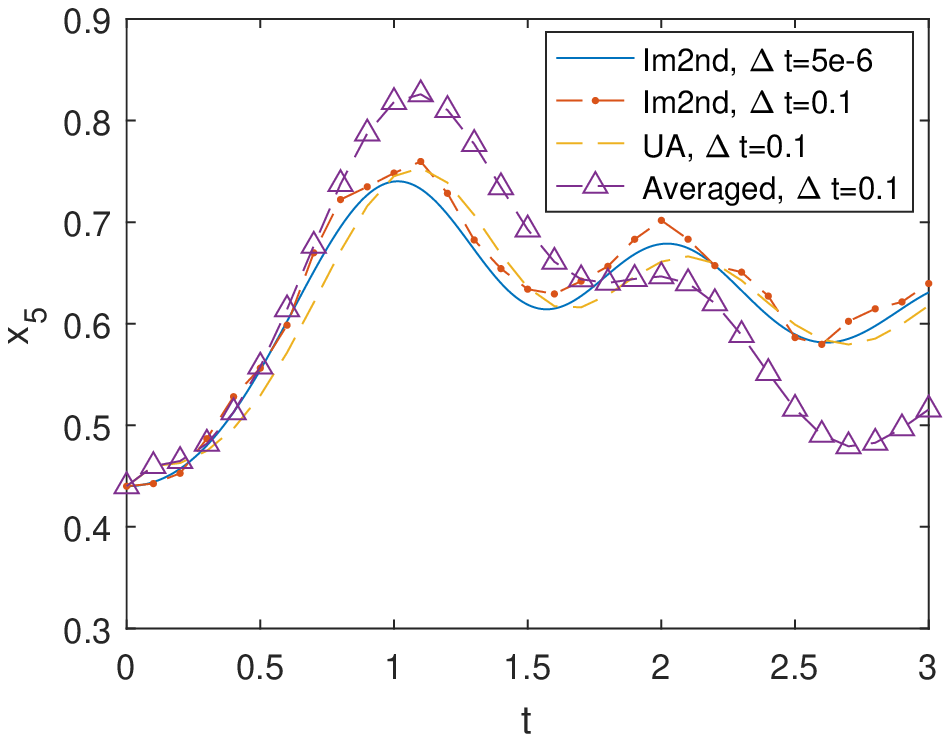}
		\label{fig:x5_avg}
	}
	\caption{$w_1$ and $w_5$ obtained by using different methods.}
	\label{fig:x1x5_avg}
\end{figure}

\textbf{Difference between $F$ and $\bar{f}^1$.} 
Simple calculations show that the right hand side of Eq.\eqref{eqn:non-stiff-nscale-1D}, i.e. 
$F(t,x)$ is approximately equal to $\bar{f}^1$. One may expect that replacing the ODE system \eqref{eqn:non-stiff-nscale-1D} by $\dot{\tilde{y}}=\bar{f}^1$ will generate a solution
that is close to the original one. We aim to investigate this issue in Fig.\ref{fig:hh4hamiltonian1e-6}.  Here, we calculate $\bar{f}^1$ and directly solve $\dot{\tilde{y}}=\bar{f}^1$, which will be referred as the averaged method.   
In Fig.\ref{fig:x1_avg} and Fig.\ref{fig:x5_avg}, we plot the first and fifth component of $\tilde{y}$, i.e., $w_1(t)$ and $w_5(t)$ obtained by several different methods. 

Fig.\ref{fig:x1_avg} shows that $w_1(t)$ is nearly a constant. Thus, the averaged method still performs well. Fig.\ref{fig:x5_avg} shows that $w_5(t)$ has oscillations. In this case, the averaged method cannot capture the right behavior of the solution, while our method can. 
We find that (1) the direct \emph{Im2nd} method with a coarse time step gives wrong results; (2) the solutions for component $w_1(t)$ obtained by the averaged method and our method agree with the reference solution; (3) the solution for component $w_5(t)$ obtained by the averaged method has large errors, while the solution for component $w_5(t)$ obtained by our method still approximates the reference solution well. This experiment shows 
that $F(t,x)$ in \eqref{eqn:non-stiff-nscale-1D} indeed captures the correct dynamics of the 
original multiscale problem, while the direct average term $\bar{f}^1$ cannot. 
\subsection{An ODE with a complicated right-hand side}
\noindent 
In the previous numerical experiments, all the integrals in our numerical schemes can be pre-computed analytically. Here, we consider a three-scale ODE, which is defined by  
\begin{equation}
\dot{x}=(1.5-\exp(\sin\frac{2\pi t}{\epsilon_{1}}+\sin\frac{2\pi t}{\epsilon_{2}}))x.
\label{ODE-3scales-ex2}
\end{equation}
In this example,  the right-hand side of the ODE \eqref{ODE-3scales-ex2} does not have analytic expressions. Thus, we use numerical quadrature rules to compute the integrals in our numerical schemes. 

Since the integrands of $f^k$, $k=1,2$ are smooth along $x$ direction, we use $8$ points in the quadrature rules to compute the integration for $\bar{f}^k$, e.g. Eqns.\eqref{mean-fn}, \eqref{mean-fn1}, \eqref{mean-fnk}, and \eqref{mean-f1}. To compute derivatives of $g$, we directly use a central difference scheme with $\Delta x=10^{-3}$. We choose the initial value $x_0=0.48$ and $(\epsilon_1,\epsilon_2)=[7\times10^{-2},11\times10^{-5}]$ in the ODE \eqref{ODE-3scales-ex2}. 
  
\begin{figure}[htbp]
	\centering
	\subfigure[]{
		\includegraphics[width=0.49\linewidth]{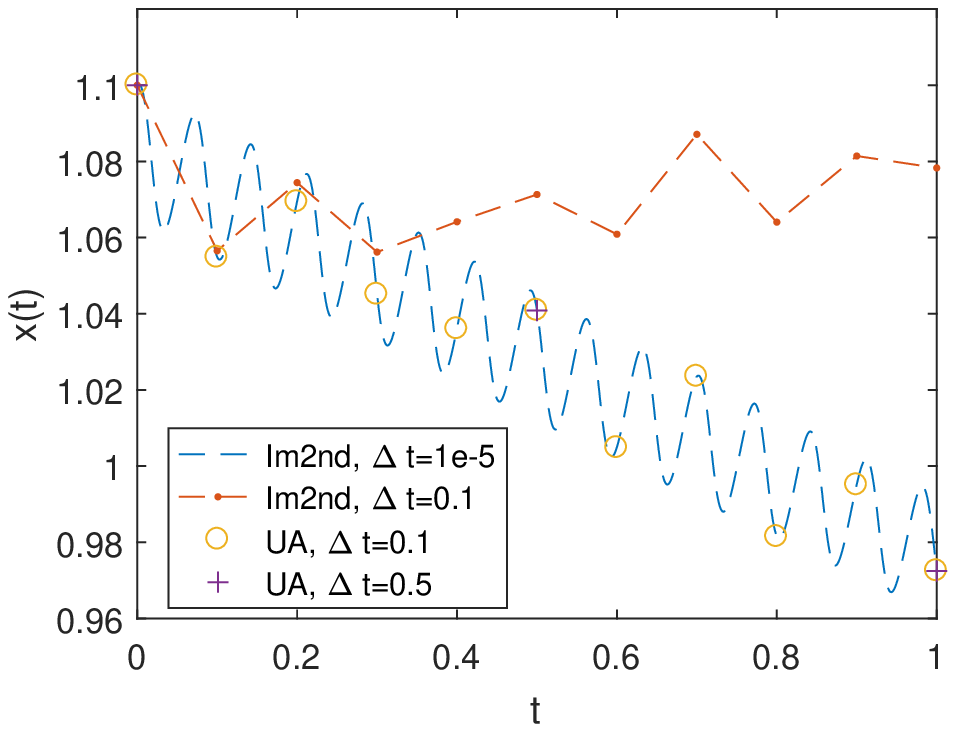}
		
		\label{fig:EgExpSin_X}
	}%
	\subfigure[]{
		\includegraphics[width=0.49\linewidth]{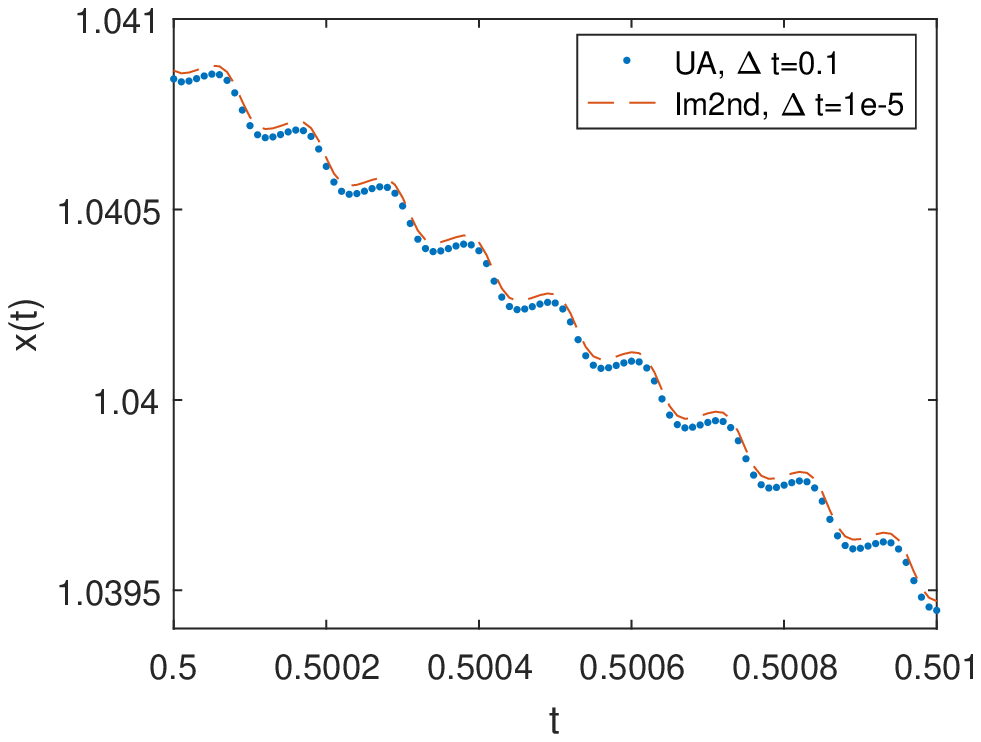}
		\label{fig:EgExpSin_Xzoom}
	}
	\caption{An ODE with a complicated right-hand side}
\end{figure}

In Fig.\ref{fig:EgExpSin_X}, we show the numerical results obtained by different methods with different time steps. We find that: (1) the direct \emph{Im2nd} method with the coarse time step $\Delta t= 0.1$ gives totally wrong results; (2) our method with very coarse time steps, $\Delta t= 0.1$ and $\Delta t= 0.5$, gives an accurate result that is comparable to that using the direct \emph{Im2nd} method with a very fine time step  $\Delta t= 10^{-5}$. This comparison again confirms that our method has a uniform accuracy in computing ODE system with multiple time-scales. In this experiment, our method with $\Delta t= 0.5$ costs $1.39$s, while the direct \emph{Im2nd} method with $\Delta t= 10^{-5}$ costs $5.97s$. More savings can be achieved if the ODE \eqref{ODE-3scales-ex2}  is described by smaller multiscale parameters $\epsilon_1$ and $\epsilon_2$. 
  
  As a byproduct of our robust numerical method (namely we can solve the complicated ODE with very large time step), we can use the representation in Eq.\eqref{xmapy} to recover the solution 
  of the ODE \eqref{ODE-3scales-ex2} in a neighborhood of the numerical solution points. 
  In Fig.\ref{fig:EgExpSin_Xzoom}, we show the recovered solution in the time domain $[0.5,0.501]$ 
  based on the solution of our method at $t=0.5$. We can see that the recovered solution agrees with 
  the reference solution and the recovered solution captures the high oscillate structure in this neighborhood. 
   
\section{Conclusions} \label{sec:Conclusion}
\noindent 
In this paper, we have successfully developed a class of robust numerical methods to solve dynamical systems with multiple time scales. These problems are difficult to solve when the multiscale parameters are small. The essential idea of our method is to represent the solution of the dynamical systems as a transformation of a slowly varying solution. Based on the scale separation assumption, we provide an efficient way to construct the transformation map and derive the dynamic equation for the slowly varying solution. Under some mild assumptions, we obtain the convergence of the proposed method. Finally, we present several numerical examples, including ODE system with three and four separated time scales to demonstrate the accuracy and efficiency of the proposed method. Numerical results show that: (1) our method is robust and accurate in solving ODE systems with multiple time scale, where the time step does not depend on the multiscale parameters; and (2) the construction of the cumulative composition maps (which deals with the multiscale information in a dimension-by-dimension fashion) is necessary while a simple average treatment leads to wrong results.  

There are two lines of work that deserve further explorations in the near future. Firstly, we shall consider to extend our proposed method to solve dynamical systems without scale separation. The idea mentioned in Remark \ref{dealwithnonscale} is a good starting point. However, we need to design fast solvers. Secondly, we are interested in extending our proposed method to solve elliptic PDEs with multiscale parameters.

\section*{Acknowledgements}
\noindent
The research of T. Hou is partially supported by the NSF Grants DMS-1613861, DMS-1907977, and DMS-1912654. The research Z. Wang is partially supported by the Hong Kong PhD Fellowship Scheme. The research of Z. Zhang is supported by Hong Kong RGC grants (Projects 27300616, 17300817, and 17300318), National Natural Science Foundation of China via grant 11601457, Seed Funding Programme for Basic Research (HKU), and Basic Research Programme (JCYJ20180307151603959) of The Science, Technology and Innovation Commission of Shenzhen Municipality. The computations were performed using the HKU ITS research computing facilities that are supported in part by the Hong Kong UGC Special Equipment Grant (SEG HKU09). 

\section*{}
\bibliographystyle{siam}

\begin{thebibliography}{10}

\bibitem{abdulle2012heterogeneous}
{\sc A.~Abdulle, E.~Weinan, B.~Engquist, and E.~Vanden-Eijnden}, {\em The
  heterogeneous multiscale method}, Acta Numerica, 21 (2012), pp.~1--87.

\bibitem{bartlett2018representing}
{\sc S.~Bartlett, P.and~Evans and P.~Long}, {\em Representing smooth functions
  as compositions of near-identity functions with implications for deep network
  optimization}, arXiv:1804.05012,  (2018).

\bibitem{chartier2015uniformly}
{\sc P.~Chartier, N.~Crouseilles, M.~Lemou, and F.~M{\'e}hats}, {\em Uniformly
  accurate numerical schemes for highly oscillatory {K}lein--{G}ordon and
  nonlinear {S}chr{\"o}dinger equations}, Numerische Mathematik, 129 (2015),
  pp.~211--250.

\bibitem{chartier2017new}
{\sc P.~Chartier, M.~Lemou, F.~M{\'e}hats, and G.~Vilmart}, {\em A new class of
  uniformly accurate numerical schemes for highly oscillatory evolution
  equations}, Found. Comput. Math., 19 (2019), pp.~1--33.

\bibitem{EfendievHou:09}
{\sc Y.~Efendiev and T.~Y. Hou}, {\em Multiscale finite element methods.
  {T}heory and applications}, Springer-Verlag, New York, 2009.

\bibitem{engquist1989particle}
{\sc B.~Engquist and T.~Hou}, {\em Particle method approximation of oscillatory
  solutions to hyperbolic differential equations}, SIAM journal on numerical
  analysis, 26 (1989), pp.~289--319.

\bibitem{ferrer1998henon}
{\sc S.~Ferrer, M.~Lara, J.~Palacian, J.~Juan, A.~Viartola, and P.~Yanguas},
  {\em The h{\'e}non and heiles problem in three dimensions. ii. relative
  equilibria and bifurcations in the reduced system}, International Journal of
  Bifurcation and Chaos, 8 (1998), pp.~1215--1229.

\bibitem{hairer2006geometric}
{\sc E.~Hairer, C.~Lubich, and G.~Wanner}, {\em Geometric numerical
  integration: structure-preserving algorithms for ordinary differential
  equations}, vol.~31, Springer Science \& Business Media, 2006.

\bibitem{henon1964applicability}
{\sc M.~H{\'e}non and C.~Heiles}, {\em The applicability of the third integral
  of motion: some numerical experiments}, The Astronomical Journal, 69 (1964),
  p.~73.

\bibitem{hou1988homogenization}
{\sc T.~Hou}, {\em Homogenization for semilinear hyperbolic systems with
  oscillatory data}, Communications on pure and applied mathematics, 41 (1988),
  pp.~471--495.

\bibitem{hou2008multiscale}
{\sc T.~Hou, D.~Yang, and H.~Ran}, {\em Multiscale analysis and computation for
  the three-dimensional incompressible {N}avier--{S}tokes equations},
  Multiscale Model. Simul., 6 (2008), pp.~1317--1346.

\bibitem{iserles2009first}
{\sc A.~Iserles}, {\em A first course in the numerical analysis of differential
  equations}, no.~44, Cambridge university press, 2009.

\bibitem{jones2012multiple}
{\sc C.~Jones and A.~Khibnik}, {\em Multiple-time-scale dynamical systems},
  vol.~122, Springer Science \& Business Media, 2012.

\bibitem{kevrekidis2003equation}
{\sc I.~Kevrekidis, C.~Gear, J.~Hyman, P.~Kevrekidid, O.~Runborg, and
  C.~Theodoropoulos}, {\em Equation-free, coarse-grained multiscale
  computation: {E}nabling mocroscopic simulators to perform system-level
  analysis}, Commum. Math Sci., 1 (2003), pp.~715--762.

\bibitem{kuehn2015multiple}
{\sc C.~Kuehn}, {\em Multiple time scale dynamics}, vol.~191, Springer, 2015.

\bibitem{perko2013differential}
{\sc L.~Perko}, {\em Differential equations and dynamical systems}, vol.~7,
  Springer Science \& Business Media, 2013.

\bibitem{shen2019nonlinear}
{\sc Z.~Shen, H.~Yang, and S.~Zhang}, {\em Nonlinear approximation via
  compositions}, arXiv:1902.10170,  (2019).

\bibitem{tao2010nonintrusive}
{\sc M.~Tao, H.~Owhadi, and J.~Marsden}, {\em Nonintrusive and structure
  preserving multiscale integration of stiff {ODE}s, {SDE}s, and {H}amiltonian
  systems with hidden slow dynamics via flow averaging}, Multiscale Model.
  Simul., 8 (2010), pp.~1269--1324.

\bibitem{tuckerman1992reversible}
{\sc M.~Tuckerman, B.~Berne, and G.~Martyna}, {\em Reversible multiple time
  scale molecular dynamics}, J. Chem. Phys., 97 (1992), pp.~1990--2001.

\end{thebibliography}

\end{document}